\documentclass{preprint}
\usepackage{amsmath}
\usepackage{amsthm}
\usepackage{amssymb}
\usepackage{graphicx}
\usepackage{mhequ}
\usepackage{mhfig}
\usepackage{times}
\usepackage{Symbols}
\usepackage{verbatim}

\newtheorem{lemma}{Lemma}[section]
\newtheorem{remark}[lemma]{Remark}
\newtheorem{theorem}[lemma]{Theorem}
\newtheorem{corollary}[lemma]{Corollary}
\newtheorem{prop}[lemma]{Proposition}
\newtheorem{proposition}[lemma]{Proposition}
\newcommand{\E}{{\mathbb E}}

\newcommand{\T}{{\mathbf T}}
\newcommand{\D}{{\mathcal D  }}

\newcommand{\cA}{\mathcal A}
\newcommand{\Lou}{L_{\mathrm{OU}}}

\newcommand{\tdif}{\tau_{\mathrm{diff}}}

\newcommand{\sscal}[1]{\langle\!\langle #1 \rangle\!\rangle}
\newcommand{\sym}{{\hbox{\scriptsize sym}}}
\def\com{{\hat C}}

\begin{document}
\title{From ballistic to diffusive behavior in periodic potentials}
\author{M. Hairer$^1$ \& G.A. Pavliotis$^2$}
\institute{Mathematics Institute, The University of Warwick, Coventry CV4 7AL, UK
\and Department of Mathematics,   Imperial College, London SW7 2AZ, UK}
\maketitle

\begin{abstract}
The long-time/large-scale, small-friction asymptotic for the one dimensional Langevin
equation with a periodic potential is studied in this paper. It is shown that the
Freidlin-Wentzell and central limit theorem  (homogenization) limits commute. We prove
that, in the combined small friction, long-time/large-scale limit the particle position
converges weakly to a Brownian motion with a singular diffusion coefficient which we
compute explicitly. We show that the same result is valid for a whole one parameter
family of space/time rescalings. The proofs of our main results are based on some novel
estimates on the resolvent of a hypoelliptic operator.
\end{abstract}
%
%
%
%
\section{Introduction and Main Results}
\label{sec:intro}
Random perturbations of dynamical systems has been the subject of intense study over the
last several decades~\cite{FW}. One of the most extensively studied randomly perturbed
dynamical systems is given by the Langevin equation modelling the interaction of
a classical particle with a heat bath at inverse temperature $\beta$:
\begin{equation}\label{e:lang}
\ddot{q} = - \nabla V(q) - \gamma \dot{q} + \sqrt{2 \gamma \beta^{-1} } \xi(t)\;.
\end{equation}
Here, $V(q)$ denotes a smooth potential, $\gamma$ is a friction coefficient which should
be interpreted as the strength of the coupling to the heat bath, and $\xi(t)$ denotes
standard $d$--dimensional white noise, i.e. a mean zero generalized Gaussian process with
correlation structure
$$
\langle \xi_i(t) \xi_j(s) \rangle = \delta_{ij} \delta(t-s), \quad i,j = 1, \dots d.
$$
There are various applications of this model to solid state physics, e.g. surface
diffusion, Josephson junctions and superionic conductors. As a result,
equation~\eqref{e:lang} has been one of the most popular stochastic models in the physics
and the mathematics literature. See, e.g., \cite{Ris84, reimann, HanTalkBork90,HelNie05HES}
and the references therein.

Various asymptotic limits for the Langevin equation~\eqref{e:lang} have been studied,
both in finite~\cite{nelson, freidlin8} and in infinite dimensions~\cite{CerFreid06a,
PavSt05a}. It is well known, for example, that for large values of the friction
coefficient $\gamma$, solutions the rescaled process
\begin{equation}\label{e:smoluch_scal}
q_{\gamma}(t)  = q(t/\gamma)
\end{equation}
converges to the solution of the Smoluchowski equation
\begin{equation}\label{e:smoluch_intro}
\gamma \dot{z} = - \nabla V(z) + \sqrt{ 2 \gamma \beta^{-1}  } \xi(t).
\end{equation}
This is usually called the {\it Kramers to Smoluchowski} limit.

Clearly, in the limit as the friction coefficient converges to zero, and for fixed finite
time intervals, we retrieve the deterministic dynamics which is governed by the
Hamiltonian system
$$
\ddot{q} = - \nabla V(q).
$$
The small $\gamma$, large-time asymptotic is much more interesting and was originally studied
by Freidlin and Wentzell~\cite{FW, FW93}. It was shown in these references that, for $d =
1$, and under appropriate assumptions on the potential, the Hamiltonian of the rescaled
process
\begin{equation}\label{e:resc_posit}
q^{\gamma} = \gamma q(t/\gamma),
\end{equation}
converges weakly, in the limit as $\gamma \rightarrow 0$, to a diffusion process on a
graph. This result was obtained for one dimensional Langevin equations with periodic
potentials--the problem we study in this paper--in~\cite{FW99}. From this limit theorem
one can infer the limiting behavior of the rescaled particle position, which actually
converges to a non-Markovian process; see Corollary~\ref{cor:integrate} and
Remark~\ref{rem:FW_q} in this paper. Results similar to those of the Freidlin-Wentzell
theory were obtained in~\cite{Sowers03, Sowers05} using singular perturbation theory.

On the other hand, when the potential is either periodic or random, and for fixed $\gamma
> 0$, the long time behavior of solutions to~\eqref{e:lang} is described by an effective Brownian
motion. Indeed, the rescaled particle position
\begin{equation}\label{e:resc_CLT}
q^{\eps}(t):= \eps q(t/\eps^2)
\end{equation}
converges weakly, in the limit as $\eps \rightarrow 0$, to a Brownian motion with a
nonnegative diffusion coefficient $D_{\gamma}$. An expression for the diffusion
coefficient can be obtained implicitly via the solution of a suitable Poisson
equation~\cite{rodenh, HairPavl04, papan_varadhan,Oll94,Kozlov89}. See also
Section~\ref{sec:homog} below.

The above limit theorem for the rescaled process $q^{\eps}(t)$ does not provide us with a
complete understanding of the long time asymptotic behavior of~\eqref{e:lang} for two
reasons. First, it does not contain any information on the time needed for the process
$q(t)$ to reach the asymptotic diffusive regime, the {\it diffusive time scale} $\tdif$.
Second, it does not provide us with any information on the dependence of the effective
diffusion coefficient $D_\gamma$ on the friction coefficient $\gamma$ and on the inverse
temperature $\beta$. The large-$\gamma$/large-$\beta$ regime is the most interesting one
from the point of view of applications and it has been studied quite extensively by means
of formal asymptotics and numerical experiments, see~\cite{sancho_al04a, sancho_al04b}
and the references therein. An asymptotic formula for the diffusion coefficient which is
valid at small temperatures was obtained rigorously by Kozlov in~\cite{Kozlov89}. The
formula obtained in that paper, however, is not valid uniformly in $\gamma$, but only for
large or intermediate values of the friction coefficient. The purpose of this paper is to
study the dependence of the diffusive time scale $\tdif$ and of the effective diffusion
coefficient $D_\gamma$ on the friction coefficient, in particular in the limit as
$\gamma$ tends to $0$, and to obtain results which are uniform in $\beta$. We also derive
various results related to the large $\gamma$ asymptotic.

\begin{figure}
\centerline{
\begin{tabular}{c@{\hspace{3pc}}c}
\includegraphics[width=2.4in, height = 2.4in]{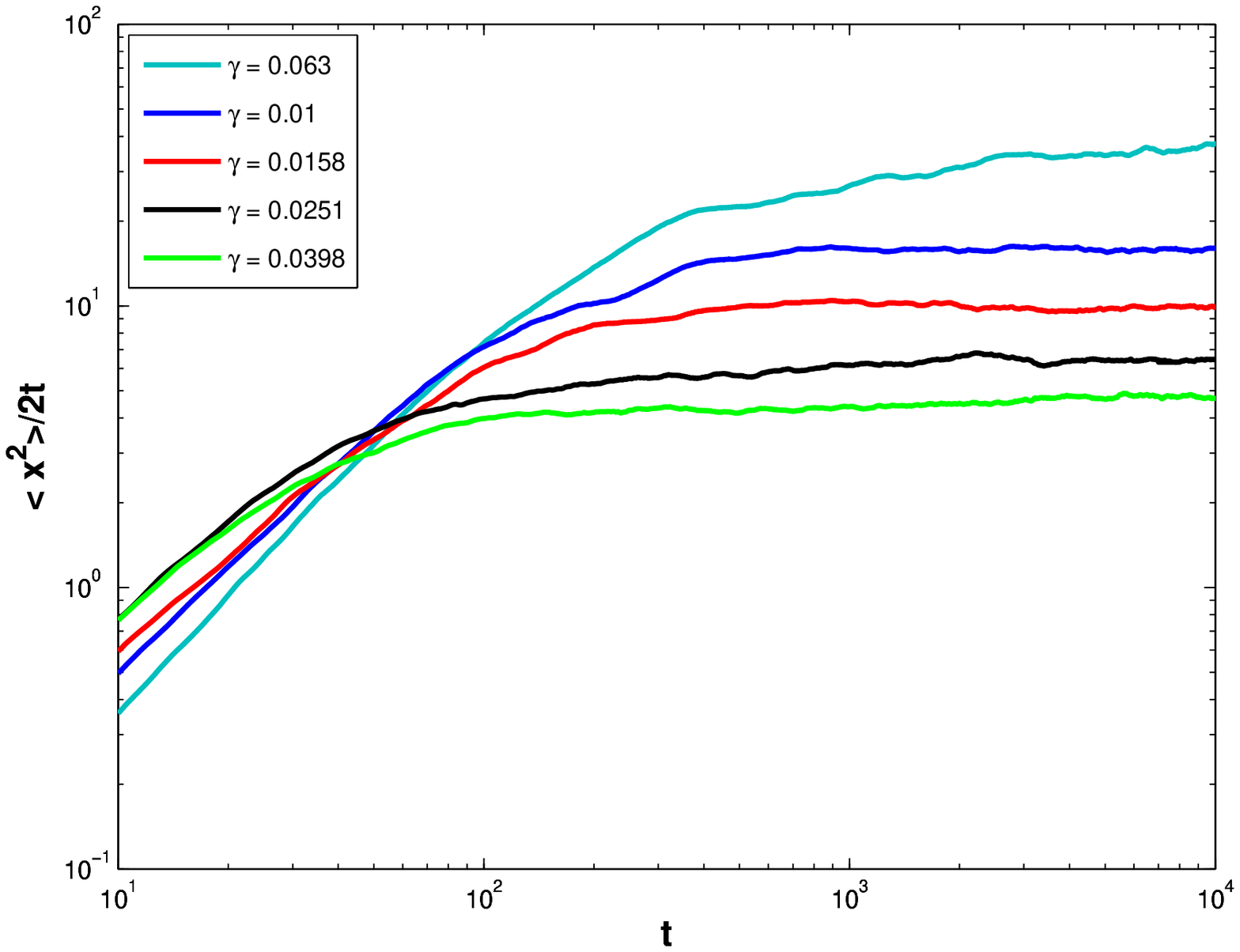} &
\includegraphics[width=2.4in, height = 2.4in]{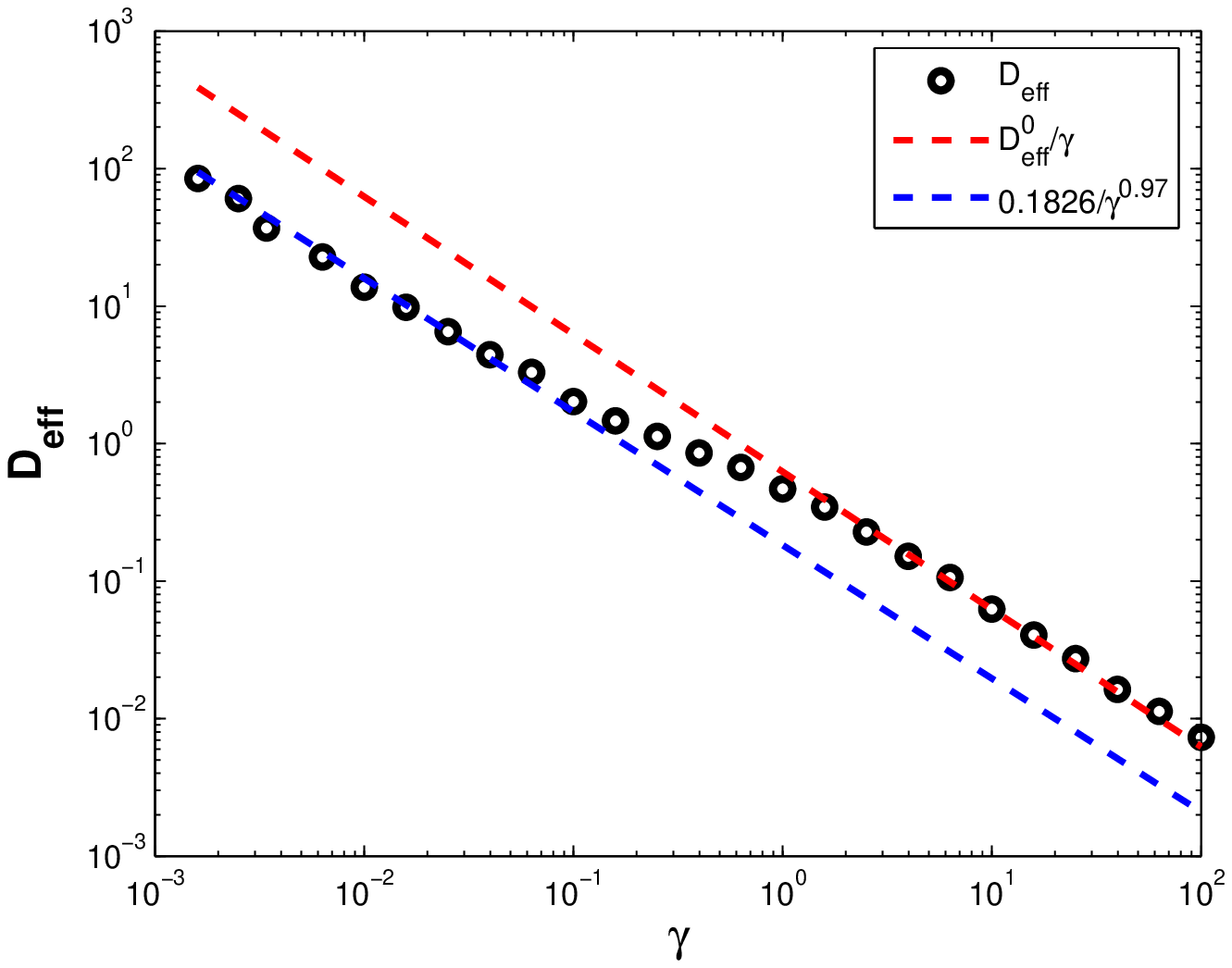} \\
a.~~  $\langle ( q(t) - \langle q(t) \rangle)^2 \rangle/2t$ vs $t$  & b.~~ $\D$ vs
$\gamma$
\end{tabular}}
\begin{center}
\caption{Second moment and effective diffusivity for various values of $\gamma$ .}
\label{fig:deff}
  \end{center}
\end{figure}

To get some intuition on the dependence of $\tdif$ and $D_\gamma$ on $\gamma$, we
calculate numerically $D_\gamma$ for the nonlinear pendulum with dissipation and noise
through the formula
$$
D_\gamma = \lim_{t \rightarrow \infty} \frac{\langle (q(t) - \langle q(t) \rangle)^2
\rangle}{2t},
$$
where $\langle \cdot \rangle$ denotes ensemble average. In Figure~\ref{fig:deff}a we plot
the second moment of the particle position divided by $2 t$ as a function of time, for
various values of the friction coefficient. In Figure~\ref{fig:deff}b we plot the
diffusion coefficient as a function of $\gamma$. All simulations were performed at a
fixed temperature $\beta^{-1} = 0.1$. The numerical simulations suggest that
\begin{equs}
\tdif \sim \frac{1}{\gamma}, \qquad &\mbox{for $\gamma \ll 1$,} \label{e:tdiff} \\
\intertext{and that} D_\gamma \sim \frac{1}{\gamma}, \qquad &\mbox{for both $\gamma \ll
1$ and $\gamma \gg 1$.}\label{e:d_gamma}
\end{equs}
The central result of this article is a rigorous justification of the above two (actually
three) scaling limits, and the explicit calculation of the prefactors for both the large
and the small $\gamma$ asymptotics of $D_\gamma$.
We will restrict our attention to the one-dimensional case.
We study the long time/small $\gamma$ asymptotic of the one-dimensional Langevin
equation
\begin{equation}\label{e:lang1}
\ddot{q} = - \partial_q V(q) - \gamma \dot{q} + \sqrt{2 \gamma \beta^{-1} } \xi(t)
\end{equation}
when $V(q)$ is a smooth, \textit{periodic} potential and $\xi(t)$ is white noise. Our first result can be
summarized in the following.
\begin{theorem}\label{thm:main}
The Freidlin--Wentzell scaling limit \eref{e:resc_posit} and the diffusive scaling limit
\eref{e:resc_CLT} `commute'. In particular, the rescaled process
$$
\eps \gamma q(t/(\gamma \eps^{2}))
$$
converges weakly, both in the $\lim_{\eps \rightarrow 0} \lim_{\gamma \rightarrow 0}$ limit
and the $\lim_{\gamma \rightarrow 0} \lim_{\eps \rightarrow 0}$ limit,
to a Brownian motion with diffusion coefficient $D^*$ given by
formula~\eqref{e:deff_FW_1} below.

Furthermore, the Kramers to Smoluchowski scaling limit~\eqref{e:smoluch_scal} and the
diffusive scaling limit~\eqref{e:resc_CLT} also 'commute': the rescaled process
$$
\eps q(t/(\gamma \eps^2))
$$
converges weakly, both in the $\lim_{\eps \rightarrow 0} \lim_{\gamma \rightarrow
\infty}$ limit and the $\lim_{\gamma \rightarrow \infty} \lim_{\eps \rightarrow 0}$
limit, to a Brownian motion with diffusion coefficient $\bar D$ given by
formula~\eqref{e:deff_smol} below.
\end{theorem}

The above theorem justifies rigorously the expressions~\eref{e:tdiff}
and~\eref{e:d_gamma}, for $\gamma \ll 1$ and $\gamma \gg 1$.

Clearly, the above theorem implies that
$$
\lim_{\gamma \rightarrow 0} \gamma D_\gamma = D^* \quad \mbox{and} \quad \lim_{\gamma
\rightarrow \infty} \gamma D_\gamma = \bar D.
$$
In fact, we can say slightly more: in Section~\ref{sec:deff_estim}, we prove the
two-sided bound
$$
\frac{D^*}{\gamma} \leq D_\gamma \leq \frac{\bar D}{\gamma} \;, \qquad\forall\, \gamma \in (0,
\infty)\;.
$$
We also compute the next order correction in the large $\gamma$ expansion of the
diffusion coefficient $D_\gamma$.

More generally, we are going to study the small-$\gamma$ asymptotic of the rescaled process
\begin{equation}\label{e:q_resc}
q^\gamma (t) = \lambda_\gamma q(t/\mu_\gamma)
\end{equation}
for a suitable one-parameter family of space-time rescalings $\lambda_\gamma, \,
\mu_\gamma$. It turns out that the ``right'' scalings -- the ones giving rise to a
non-trivial limiting process -- are of the form
\begin{equation}\label{e:resc_alpha}
\lambda_\gamma = \gamma^{1 + \alpha}, \quad \mu_\gamma = \gamma^{1 + 2 \alpha}, \quad
\alpha \in [0, \infty).
\end{equation}
Note that the case $\alpha = 0$ corresponds to the Freidlin-Wentzell rescaling
\eref{e:resc_posit}, whereas the case $\alpha = \infty$ corresponds to the diffusive
rescaling \eref{e:resc_CLT}. Our second result is the following.
\begin{theorem}\label{thm:main_2}
Assume that the Markov process $(q(t), \, p(t))$ is stationary on $\T \times \R$. Then
the rescaled process $q^{\gamma}(t)$ defined in~\eqref{e:q_resc} converges weakly to a
Brownian motion for every $\alpha \in (1/2, +\infty)$. The diffusion coefficient
coefficient of the limiting Brownian motion is independent of $\alpha$ and is given
by~\eqref{e:deff_FW_1}.
\end{theorem}
\begin{remark}
We believe that this is the theorem is also true for $\alpha \in (0, 1/2)$. However, we
haven't been able to prove this. See also Remark~\ref{rem:l4_estim} below.
\end{remark}
\begin{remark}
The stationarity assumption is not necessary and can be replaced with the assumption that
the distribution of the initial condition has an $L^2$ density with respect to the
Maxwell-Boltzmann distribution $\mu(dp \, dq) = Z^{-1} \exp(- \beta H(p,q)) \, dp \, dq$.
For purely technical reasons it seems to be more difficult to obtain the same result for
deterministic initial conditions.
\end{remark}
The, perhaps, surprising result is that the diffusion coefficient is {\it independent of
the exponent} $\alpha$: as long as we are at length and time scales which are long
compared to the Freidlin-Wentzell length and time scales, the particle performs an
effective Brownian motion {\it with the same diffusion coefficient}.

\begin{remark}\label{rem:gamma_large}
A similar result holds for the large $\gamma$ limit:  Under the assumption of
stationarity, we have that
$$
\lim_{\gamma \rightarrow \infty}\gamma^{-\alpha} q(t \gamma^{1 + 2 \alpha}) = \sqrt{2
\bar D} \, W(t)
$$
weakly on $C([0,T], \R)$ for every $\alpha > 0$, where $\bar D$ is given by
formula~\ref{e:deff_smol} below.
\end{remark}

Similar scalings to the one considered in~\eqref{e:q_resc} were considered for the
passive tracer dynamics
$$\dot{q} = v(q) + \sqrt{2 \sigma} \xi, \quad \nabla \cdot v = 0 $$
by Fannjiang in~\cite{Fann01}. There it was shown that the diffusive time scale depends
crucially on the ergodic properties of the vector field $v(q)$ on $\T^d$. On the
contrary, for the problem studied in this paper, the small $\gamma$ asymptotic of
$\tdif$ and $D$ are independent of the specific properties of the potential $V(q)$.
This is because the Hamiltonian vector field can never generate an ergodic flow on the
phase space $\T \times \R$ due to the presence of the integral of the energy.

The proofs of Theorems~\ref{thm:main} and~\ref{thm:main_2} are based on a careful analysis
of the generator of the Markov process $(q(t), p(t))$ on $\T \times \R$.
It turns out to be notationally more convenient to study the the rescaled generator of
\eref{e:lang1}
\begin{equation}\label{e:generator}
L_\gamma = \frac{1}{\gamma} \cA + \Lou,
\end{equation}
on $\T \times \R$, where $\cA = p\d_q - V'(q)\d_p$ is the Liouville operator describing the
unperturbed deterministic dynamic and $\Lou =
\beta^{-1}\d_p^2 - p\d_p$ is the generator of the Ornstein-Uhlenbeck process describing the
interaction with the heat bath.

The main technical results which are needed for the proof of Theorem~\ref{thm:main} are
an estimate on the resolvent of $L_\gamma$, as well as estimates on derivatives of
solutions to Poisson equation of the form $-L_\gamma u =h$. We obtain an estimate on the
semigroup generated by $L_\gamma$ which is independent of $\gamma$:
\begin{theorem}\label{theo:sgap}
There exist constants $C$ and $\alpha$ independent of $\gamma$ such that
\begin{equ}[e:sgbound]
\|e^{L_\gamma t} f\| \le Ce^{-\alpha t} \|f\|\;,
\end{equ}
holds for every $t>0$, every $\gamma < 1$, and every  $f \in \L^2(\mu)$ such
that $\int f d \mu = 0$, where $\mu( dp \, dq) = Z^{-1} \exp(- \beta H(q,p)) \, dp \, dq$.
\end{theorem}
The Poisson equation that we need to analyze is
\begin{equation}\label{e:poisson}
- L_\gamma \phi_\gamma = p.
\end{equation}
The boundary conditions for this PDE are that the solution in periodic in $q$ and that it
belongs to $L^2(\mu)$. Our estimate on derivatives of $\phi_\gamma$ is uniform in
$\gamma$:
\begin{prop}\label{prop:estim_l4}
Assume that $V(q)$ is smooth and let $\phi_\gamma$ be the solution to~\eqref{e:poisson}.
The there exists a constant $C$ which is independent of $\gamma$ such that
\begin{equation}
\|\phi_{\gamma} \|^2 + \|\d_p \phi_{\gamma}\|^2 + \|\d_q \phi_{\gamma}\|^2 + \gamma \Big(
\|\d_p^2 \phi_{\gamma}\|^2 + \|\d_p \d_q \phi_{\gamma}\|^2 + \|\d_q^2 \phi_{\gamma}\|^2
\Big) \le C\;,
\end{equation}
independently of $\gamma$. Furthermore, $\partial_p \phi_\gamma$ is an element of
$L^4(\mu)$ and
\begin{equation}\label{e:estim_l4}
\|\partial_p \phi_\gamma \|_{L^4(\mu)} \leq C (1 + \gamma^{-1/4})\;.
\end{equation}
\end{prop}
\begin{remark}\label{rem:l4_estim}
We believe that estimate~\ref{e:estim_l4} should actually be uniform in $\gamma$. However, we
haven't been able to prove this. The reason for this is that we obtain \eref{e:estim_l4}
as a consequence of Sobolev embedding, but $\d_p^2 \phi_\gamma$ is not uniformly
bounded in any weighted $L^2$ space.
\end{remark}
The proof of estimates~\ref{e:sgbound} and~\ref{e:estim_l4} is based on the commutator
techniques that were developed recently by Villani~\cite{Vil04HPI}. Somewhat similar
estimates to the ones we prove in this paper were recently derived by H\'erau in~\cite{Her07}.

The rest of the paper is organized as follows. In Section~\ref{sec:critical} we analyze
the Freidlin-Wentzell scaling~\eqref{e:resc_posit}. In Section~\ref{sec:homog}, we then
study the diffusive scaling~\eqref{e:resc_CLT}. In Section~\ref{sec:deff_estim} we obtain
upper and lower bounds on the diffusion coefficient and we study the large $\gamma$
asymptotic. The intermediate scalings~\eqref{e:q_resc} for $\alpha \in (0, +\infty)$ are
investigated in Section~\ref{sec:inter}.  The necessary estimates on the resolvent of the
generator $L_{\gamma}$ are presented in Section~\ref{sec:resolvent}.

\subsection*{Acknowledgements} The authors would like to thank R. Sowers, E. Vanden Eijnden
and A.M. Stuart for useful discussions and comments.
%
%
%
\section[Critical Scaling: the $\alpha = 0$ case]{Critical Scaling: the $\pmb{\alpha = 0}$ case}
\label{sec:critical}

Let us rewrite the Langevin equation~\eqref{e:lang} in one space dimension as a first order system:
\minilab{e:lang_syst}
\begin{equs}
d q(t) &= p(t) \, dt\;, \label{e:lang_1} \\
d p(t) &= -  \partial_q V(q(t)) \, dt - \gamma p(t) \, dt +\sqrt{2 \gamma \beta^{-1}} \, d
W\;, \label{e:lang_2}
\end{equs}
where $V(q)$ is a smooth periodic potential with period $1$ and $W(t)$ is a standard
one-dimensional Wiener process.

This section is devoted to the study of the critical scaling
\begin{equ}[e:defxgamma]
q^\gamma(t) = \gamma q(t/\gamma)\;,
\end{equ}
which corresponds to the limiting case which is not covered by Theorem~\ref{thm:main_2}.
It turns out that under this scaling, $q^\gamma$ does \textit{not} converge to a Brownian
motion, but to a non-Markovian process that will be described in this section.

\begin{figure}
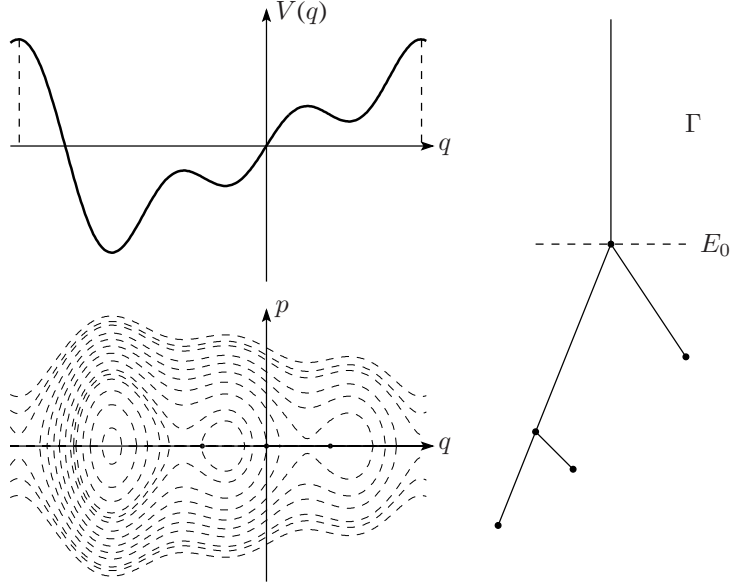

\begin{center}
\mhpastefig{Graph}
\end{center}
\caption{Example of a potential with the orbits of the Hamiltonian flow and the
corresponding Freidlin-Wentzell graph $\Gamma$.}\label{fig:FW}
\end{figure}

The behavior at the critical scaling can be understood with the help of the
Freidlin-Wentzell theory of averaging for small random perturbations of a Hamiltonian
system \cite{FW93,FW98,FW99}. Recall that one can associate to a Hamiltonian system on
the symplectic manifold $\CM = \T \times \R$ a graph $\Gamma$ in such a way that every
point in the graph corresponds to a connected component of a level set of $H$. Vertices
of the graph correspond to level sets containing a critical point of $H$. See
Figure~\ref{fig:FW} for an example.

We identify points on the graph $\Gamma$ with elements of $\R \times \Z$ by
ordering the edges of the graph
and taking the value of the Hamiltonian as a local coordinate along each edge. We denote by $\tilde H \colon \CM \to \Gamma \approx \R \times \Z$ the
`extended' Hamiltonian which associates each point to its energy, together with the
number of the edge to which the corresponding connected component belongs.

Denoting by $\lambda$ the Lebesgue measure on $\CM$, the measure $\tilde \lambda = \tilde H^* \lambda $
on $\Gamma$ then has a density with respect to Lebesgue measure on $\Gamma$, which we
denote by $T(z)$. The notation $T(z)$ is justified by the fact that it is actually equal to the
period of the orbit corresponding to
the point $z$. It is therefore a straightforward exercise to see that
\begin{equ}
T(z) \approx |\log (z - z_0)|\;,
\end{equ}
near the vicinity of a critical orbit $z_0$ which corresponds to a maximum of the potential.
For a point $z \in \Gamma$, denote by $\ell_z$ the measure on $\tilde H^{-1}(z)$
which is such that
\begin{equ}
\int_\CM f(x) \, \lambda(dx) = \int_\Gamma \int_{\tilde H^{-1}(z)} f(x)\,\ell_z(dx)\, dz\;,
\end{equ}
for every integrable function $f\colon \CM \to \R$. The measure $\ell_z$ is not a probability
 measure but has mass $T(z)$. With this notation at hand, we define the function
\begin{equ}
S(z) = \int_{\tilde H^{-1}(z)} p^2\, \ell_z(dx)\;,
\end{equ}
where we used the notation $x = (q,p)$ for elements of $\CM$. The function $S(z)$ has non-trivial limits
as $z$ approaches the vertices of $\Gamma$. Note that these limits are in general different
for different ways of approaching the same vertex, so that $S$ is discontinuous on $\Gamma$.
It is also possible to check \cite{FW}  that $S$ satisfies the relation $S'(z) = T(z)$ in the interior of the edges.

The main result of \cite{FW99} is then
\begin{theorem}\label{theo:FW}
Let $X^\gamma(t)$ be defined by $X^\gamma(t) = (p(t/\gamma), q(t/\gamma))$, where $(p,q)$
is a  solution to \eref{e:lang_syst}. Then, the process $\tilde H(X^\gamma (t))$
converges weakly to a Markov process $Y$ on $\Gamma$ whose generator is given by
the expression
\begin{equ}[e:defL]
Lv(z) = {1 \over \beta T(z)} {d \over dz} \Bigl(S(z) {dv(z) \over dz}\Bigr)
-  {S(z) \over T(z)} {dv(z)\over dz}\;,
\end{equ}
for $z$ in the interior of the edges of $\Gamma$. The domain of $L$ consists of functions
$v$ such that the above expression is square integrable, and such that at each interior vertex, the derivatives of $v$ along the edges satisfy the `gluing' conditions
\begin{equ}
\sum_{k \sim z_0} \sigma(z_0,k) \lim_{z \to_k z_0} S(z) {dv(z) \over dz} = 0.
\end{equ}
Here, $z_0$ denotes an interior vertex of $\Gamma$, the sum runs over all edges $k$ adjacent to $z_0$,
and $z \to_k z_0$ means that $z$ converges to $z_0$ along the edge $k$.
The factor $\sigma(z_0,k)$ is equal to $1$ if $H(z) > H(z_0)$ for $z$ in the $k$th edge and $-1$ otherwise.
\end{theorem}

Note that the gluing conditions are such that the process $Y$ is reversible
with respect to the probability measure $\mu_\beta(dz) = Z_\beta^{-1} e^{-\beta z} \tilde \lambda(dz)
= Z_\beta^{-1} e^{-\beta z} T(z)\,dz$ on $\Gamma$. This in turn is precisely
the push-forward under $\tilde H$ of the probability measure $Z_\beta^{-1} e^{-\beta H(x)}\, \lambda(dx)$ on $\CM$
which is invariant for the process $X^\gamma$.

\begin{corollary}\label{cor:integrate}
Let $f \colon \CM \to \R$ be smooth with at most polynomial growth and define $\bar f \colon \Gamma \to \R$ by
\begin{equ}
\bar f(z) = {1\over T(z)} \int_{\tilde H^{-1}(z)} f(x)\, \ell_z(dx)\;.
\end{equ}
Then, the process $\int_0^t f(X^\gamma(s))\,ds$ converges weakly to the process $\int_0^t \bar f(Y(s))\, ds$.
\end{corollary}

\begin{proof}
If $f \equiv 0$ in a neighborhood of the critical orbits, the result follows from a
standard averaging argument which will not be reproduced here. We refer to \cite{} for a
similar calculation. We can now construct smooth functions $f^\delta$ such that $f^\delta
\equiv 0$ in a $\delta$-neighborhood of the critical orbits and $f^\delta = f$ outside of
a $2\delta$-neighborhood of the critical orbits. The result then follows immediately from
the fact that there exists a function $h$ with $\lim_{\delta \to 0}h(\delta)=0$ such that
the expectation of the time that the process $X^\gamma$ spends in the region where
$f^\eps \neq f$ is bounded by $h(\delta)$, uniformly in $\gamma$ (see \cite[p.~294]{FW}).
\end{proof}

\begin{remark}\label{rem:FW_q}
An important particular case of Corollary~\ref{cor:integrate} is that of $f(p,q) = p$. It
shows that the process $q^\gamma$ defined in \eref{e:defxgamma} converges weakly to the
process
\begin{equ}
q^*(t) = \int_0^t \bar p(Y(s))\, ds\;,
\end{equ}
where the function $\bar p \colon \Gamma \to \R$ is defined from $p$ as in Corollary~\ref{cor:integrate}.
\end{remark}

It is clear that the process $q^*$ is not Markov by itself, but requires the computation of $Y$ first.
Note also that the function $\bar p(z)$ vanishes identically for values of $z$ corresponding to closed
orbits, so that the process $q^*$ is constant on intervals of time of positive length.
On values of $z$ for which the orbits are open, one has
\begin{equ}[e:exprp]
\bar p(z) = \pm {1 \over T(z)}\;,
\end{equ}
since the average velocity is given by the size of the torus (which was set to $1$), divided by
the period of the orbit.

Denote by $\CP_t$ the semigroup over $\Gamma$ generated by $L$.
It follows from the central limit theorem for additive functionals of reversible Markov
processes~\cite{kipnis} that the process
$\eps q^*(t/\eps ^2)$ converges weakly as $\eps \to 0$ to a Brownian motion with diffusivity
given by
\begin{equ}
D^* = \int_0^\infty \int_\Gamma \bar p(z) \CP_t \bar p(z)\, \mu_\beta(dz)\, dt\;.
\end{equ}
Since $L$ is self-adjoint in $\L^2(\Gamma, \mu_\beta)$ and has a spectral gap, this integral
converges and is given by
\begin{equ}\label{e:deff_fw}
D^* = -\scal{\bar p, L^{-1} \bar p}_\beta\;,
\end{equ}
where we denoted by $\scal{\cdot,\cdot}_\beta$ the scalar product in $\L^2(\Gamma, \mu_\beta)$.

It turns out that in our case, this expression can be computed in a very explicit way. Note that in
$\L^2(\Gamma, \mu_\beta)$, one has
$L = - A^*A$, where the first order differential operator $A$ is given by
\begin{equ}
A v (z) = \sqrt{S(z) \over \beta T(z)} {dv(z) \over dz} \equiv a(z) {dv(z)\over dz}\;.
\end{equ}
The domain of $A$ consists of all continuous functions $f$ on $\Gamma$ such that
$f$ is weakly differentiable in the interior of each edge and such that
$Af \in \L^2(\Gamma, \mu_\beta)$.

The adjoint of $A$ is given by the operator that acts in the interior of the edges of $\Gamma$ like
\begin{equs}
A^* w(z) &= - {e^{\beta z} \over T(z)} {d\over dz} \bigl(T(z) e^{-\beta z} a(z) w(z)\bigr) \\
&=- {d\over dz} \bigl(a(z) w(z)\bigr) - w(z) a(z) \Bigl({T'(z) \over T(z)} - \beta\Bigr)\;,
\end{equs}
endowed with the `boundary conditions'
\begin{equ}[e:BCA]
\sum_{k \sim z_0} \sigma(z_0,k)\lim_{z \to_k z_0} T(z) a(z) w(z) = 0\;.
\end{equ}
Here, we used the same notations as in the statement of Theorem~\ref{theo:FW}.
One then has the following variational formulation of $D^*$:
\begin{equ}[e:varD]
D^* = \inf \bigl\{\|g\|_\beta^2\,|\, A^*g = \bar p \bigr\}\;.
\end{equ}

Functions satisfying the relation $A^* g = \bar p$ are of the form
\begin{equ}
g(z) = {\sqrt \beta e^{\beta z} \over \sqrt{ S(z) T(z)}} \Bigl(V_k + \int_{z_0}^z T(z) \bar p(z) e^{-\beta z}\,dz \Bigr)\;,
\end{equ}
where we denote by $k$ the index of the edge to which $z$ belongs and by $z_0$ the vertex with the
lowest energy adjacent to that edge. The constants $V_k$ are determined
by the requirements that $g$ satisfies the conditions \eref{e:BCA} and that $g \in \L^2$.
By \eref{e:varD}, remaining degrees of freedom should be dealt with by minimising over
$\|g\|_\beta$.

In our case, the graph $\Gamma$ contains two infinite edges and a number of finite ones.
Since $\bar p$ vanishes on the finite edges and is given by \eref{e:exprp} on the two
infinite edges, it follows that the function $g$ minimizing \eref{e:varD} is given by
\begin{equ}
g(z) = \sigma(z) {1 \over  \sqrt{\beta S(z) T(z)}} \;,
\end{equ}
where the function $\sigma(z)$ vanishes on all the finite edges and is equal to $\pm 1$
on the infinite edges, with the same sign as $\bar p$. Therefore, we finally obtain for
$D^*$ the expression
\begin{equ}\label{e:deff_FW_1}
D^* = {2 \over \beta Z_\beta} \int_{E_0}^\infty {e^{-\beta z} \over S(z)}\, dz\;,
\end{equ}
where $E_0$ is the energy of the vertex at which the two infinite edges join.
(The reason for the factor $2$ in front of the above expression
is that there are exactly two infinite edges starting at $E_0$.)
The function $S(z)$ is asymptotic to $2zT(z) \approx \sqrt{2z}$ at infinity and converges to a non-zero constant
as $z \to E_0$. Furthermore, the partition function $Z_\beta$ behaves like $T_0/\beta$ for large
values of $\beta$ (here $T_0$ is the value of $T(z)$ as $z$ approaches the orbit where
the energy attains its global minimum).

In order to compute the behavior of $Z_\beta$ for small values of $\beta$, we use the
fact that $T(z) \approx {1\over \sqrt{2z}}$ for large values of $z$. Therefore
\begin{equ}
Z_\beta \approx \int_0^\infty {e^{-\beta z} \over \sqrt{2z}} \,dz = \sqrt{\pi \over 2\beta}\;.
\end{equ}
A similar calculation allows to evaluate the behavior of the integral over $e^{-\beta
z}/S(z)$ for small values of $\beta$. Collecting these asymptotic estimates, one obtains
\begin{equs}[2]
D^* &\approx {2 \over \beta} &\qquad  \beta &\to 0\;, \\
D^* &\approx {2 e^{-\beta E_0} \over \beta T_0 S(E_0)} &\qquad \beta &\to \infty\;.
\end{equs}

\begin{remark}
It is unsurprising to see that the high-temperature limit $\beta \to 0$ coincides with
the result that one obtains when $V \equiv 0$.
\end{remark}
%
%
%
%
%
\section[The Central Limit Theorem Regime: The $\alpha = \infty$ Case.]
{The Central Limit Theorem Regime: The $\pmb{\alpha = \infty}$ Case.} \label{sec:homog}
Just as in the previous section, the long time behavior of solutions to \eqref{e:lang}
for a fixed value of $\gamma$ is governed by an effective Brownian motion. Indeed, the
following central limit theorem holds \cite{rodenh, papan_varadhan, HairPavl04,
Kozlov89}.
\begin{theorem}
Let $V(q) \in C^{\infty}_{per}(\T)$ and define the rescaled process
$$
q^\eps (t) := \eps q(t/\eps^2).
$$
Then $q^{\eps} (t)$ converges weakly, on $C([0,T], \R)$, in the limit as $\eps
\rightarrow 0$, to a Brownian motion with diffusion coefficient
\begin{equation}
D_\gamma = \frac{1}{\gamma} \int_{\T \times \R} p  \phi_\gamma \mu(dp \, dq)\;,
\label{e:efdif}
\end{equation}
where $\mu (dp \, dq) = Z^{-1}\exp(-\beta H(p,q))\, dp\, dq$, and the function
$\phi_\gamma$ is the unique mean-zero solution of the Poisson equation
\begin{equation}
-L_\gamma \phi_\gamma = p\;. \label{e:cell}
\end{equation}
Here $L_\gamma$ is the rescaled generator defined in~\eqref{e:generator} and
$\phi_\gamma$ is periodic in $q$ and an element of $L^2(\mu)$.
\end{theorem}
\begin{remark}
This theorem is valid in arbitrary dimensions. It is also valid when the force field in
\eqref{e:lang} is not the gradient of a scalar function, provided that $\mu(dp \, dq)$ is
replaced by the corresponding invariant measure; see \cite{HairPavl04}.
\end{remark}
The main result of this section is that, in the limit as the friction coefficient
$\gamma$ tends to $0$, the rescaled effective diffusion coefficient given
by~\eqref{e:efdif} converges to the Freidlin-Wentzell effective
diffusivity~\eqref{e:deff_fw}.
\begin{prop}\label{cor:limit}
One has $\lim_{\gamma \to 0} \gamma D_\gamma = D^*$, where $D^*$ is obtained by
\eref{e:deff_fw}.
\end{prop}

\begin{proof}
Denote as before by $L_\gamma$ the generator of the critically rescaled dynamic
\eref{e:defxgamma} given in eqn.~\eqref{e:generator}  and by $\CP_t^\gamma$ the
corresponding semigroup acting on $\L^2(\CM, \mu)$. Denote furthermore as previously by
$L$ the generator of the limiting Feeidlin-Wentzell dynamic \eref{e:defL} and by $\CP_t$
the corresponding semigroup acting on $\L^2(\Gamma, \mu)$. Finally, we introduce the
averaging operator $\Pi$ defined (on continuous functions $f \colon \CM \to \R$) by
\begin{equ}[e:average]
(\Pi f)(z) = {1\over T(z)} \int_{\tilde H^{-1}(z)} f(y)\,\ell_z (dy)\;,\qquad z \in
\Gamma\;.
\end{equ}
Note that $\Pi f$ is a function from $\Gamma$ to $\R$. Furthermore, it is immediate that
$\Pi$ is a contraction from $\L^2(\CM,\mu)$ to $\L^2(\Gamma, \mu)$ and can therefore be
extended uniquely to all
 of $\L^2(\CM,\mu)$.

We also define the isometric embedding operator $\iota\colon \L^2(\Gamma,\mu)\to\L^2(\CM,
\mu)$ by
\begin{equ}
(\iota f)(x) = f(\tilde H(x))\;.
\end{equ}
With these notations, one has $D^* = \scal{\Pi p, L^{-1} \Pi p}$ and $\gamma D_\gamma =
\scal{p, L_\gamma^{-1} p}$, so that the result follows if one can show that the strong
limit in $\L^2(\CM,\mu)$
\begin{equ}[e:convergenceFWresolvent]
\lim_{\gamma \to 0} L_\gamma^{-1} f = \iota L^{-1} \Pi f\;,
\end{equ}
holds for every element $f \in \L^2(\CM,\mu)$ such that $\int f(x) \,\mu(dx) = 0$.

This will be the consequence of the following two lemmas:
\begin{lemma}\label{lem:conv}
For every function $f \in L^2(\CM, \mu)$, the limit $\lim_{\gamma \to 0} \CP_t^\gamma f =
\iota \CP_t \Pi f$ holds in $\L^2(\CM, \mu)$.
\end{lemma}
\begin{proof}
Assume first that $f$ is bounded and continuous. It then follows from
Corollary~\ref{cor:integrate} that $\lim_{\gamma \to 0} \bigl(\CP_t^\gamma f\bigr)(x) =
\bigl(\iota \CP_t \Pi f\bigr)(x)$ for every $x \in \CM$. The claim then follows from
Lebesgue's dominated convergence theorem. The fact that the claim holds for every $f \in
L^2(\CM, \mu)$ is now a simple consequence of the density of bounded continuous
functions, together with the fact that $\CP_t^\gamma$ is a contraction operator in
$L^2(\CM, \mu)$.
\end{proof}

This, together with Theorem~\ref{theo:sgap} yields:
\begin{lemma}\label{lem:experglim}
There exist constants $C$ and $\alpha$ (independent of $\gamma < 1$) such that
\begin{equ}
\|\CP_t^\gamma f\| + \|\CP_t f\| \le C e^{-\alpha t} \|f\|\;,
\end{equ}
for every $f \in \L^2(\Gamma, \mu)$ such that $\int f(x)\,\mu(dx) = 0$, and for every
$t>0$.
\end{lemma}

\begin{proof}
The bound on $\|\CP_t^\gamma f\|$ is precisely the one given in Theorem~\ref{theo:sgap}.
Since this bound is uniform in $\gamma$, the bound on $\CP_t f$ follows at once from
Lemma~\ref{lem:conv}.
\end{proof}

We now have all the necessary ingredients for the proof of
\eref{e:convergenceFWresolvent}. Fix $\eps > 0$ and choose $T$ sufficiently large such
that
\begin{equ}
\Bigl\|L_\gamma^{-1} f - \int_0^T \CP_t^\gamma f\,dt \Bigr\| \le \eps\;,\qquad
\Bigl\|L^{-1} \Pi f - \int_0^T \CP_t \Pi f\,dt \Bigr\| \le \eps\;.
\end{equ}
Such a $T$ can be chosen independently of $\gamma$ by Lemma~\ref{lem:experglim}. On the
other hand, it follows from Lemma~\ref{lem:conv} and Lebesgue's dominated convergence
theorem that
\begin{equ}
\lim_{\gamma \to 0} \Bigl\|\int_0^T\Bigl(\CP_t^\gamma f - \iota \CP_t \Pi f\Bigr)\,dt
\Bigr\| = 0\;,
\end{equ}
and the result follows.
\end{proof}
\begin{remark}
Following the methodology advertised in \cite{PavlSt06b}, it would be satisfying to
obtain an expansion for $\phi_\gamma$ of the type $\phi_\gamma = \phi_0 + \gamma \phi_1 +
\rho$ for some error term $\rho$ and therefore to get better explicit control over the
convergence in \eref{e:convergenceFWresolvent}. The problem with this approach is the
loss of regularizing properties of the resolvent $L_\gamma^{-1}$ as $\gamma \to 0$. In
particular, the limiting function $\phi_0$ is not $\CC^\infty$, but only Lipschitz
continuous. As a consequence, the first corrector is not of order $\gamma$, but expected
to be of order $\gamma^{1/2}$, see \cite{Sowers03, Sowers05}, thus leading to a breakdown
of the naive perturbative expansion.
\end{remark}
%
%
%
%
%
\section{Estimates on the Effective Diffusion Coefficient}
\label{sec:deff_estim}
In this section we present some estimates on the diffusion coefficient $D_\gamma$ defined
in~\eqref{e:efdif}. To state the upper bound we need to define the diffusion coefficient
for the Smoluchowski equation
\begin{equation}
\dot{z} = -\partial_z V(z) + \sqrt{\frac{2}{\beta}} \, \xi(t),
\label{e:smoluch1}
\end{equation}
with $V(z)$ being the smooth periodic potential in \eqref{e:lang}. It is well known, see e.g.
\cite{Oll94} or \cite[Ch. 13]{PavlSt06b}, that the rescaled process $\eps z(t/\eps^2)$
converges weakly in the limit as $\eps \rightarrow 0$ to $\sqrt{2 \bar D} W(t)$ where $W(t)$
is a standard Brownian motion and the diffusion coefficient is given by the formula
\begin{equation}\label{e:deff_smol}
\bar D = \beta^{-1} \int_{\T} |1 + \partial_q \chi|^2 \, \nu (dq) =: \beta^{-1} \|1 +
\partial_q \chi \|^2,
\end{equation}
where
$$
\nu (d q) = \frac{1}{Z} e^{- \beta V(q)} \, dq, \quad Z = \int_{\T} e^{- \beta V(q)} \,
dq\;,
$$
and the function $\chi$ is the solution to the Poisson equation
\begin{equation}\label{e:cell_smol}
\bar L \chi = \partial_q V(q), \quad \bar L = - \partial_q V(q) \partial_q + \beta^{-1}
\partial^2_q\;,
\end{equation}
equipped with periodic boundary conditions. It is well known that $\bar D \leq \beta^{-1}$. The upper
bound in the theorem below shows that diffusion for the Langevin dynamics is depleted even
further.
\begin{prop}
Let $D^*$ be as in \eref{e:deff_FW_1} and let $\bar D$ be as above. Then, the bound
\begin{equation}
\frac{D^*}{\gamma} \leq D_\gamma  \leq \frac{\bar D}{\gamma} \;,
\label{e:deff_a_1}
\end{equation}
is valid for every $\gamma \in (0,\infty)$.\label{thm:upper_bd}
\end{prop}
\begin{proof}
We multiply equation~\eqref{e:cell} by a smooth test function $\psi \in L^2(\mu)$
to obtain
\begin{equation}\label{e:weak}
\frac{1}{\gamma}\int_{\T \times \R} \phi_\gamma \cA \psi \, \mu(dp \, dq)  +  \beta^{-1}
\int_{\T \times \R} \partial_p \phi_\gamma \partial_p \psi \, \mu(dp \, dq) = \int_{\T
\times \R} p \psi \, \mu(dp \, dq).
\end{equation}
We choose a test function which is independent of $p$, $\psi = \psi(q)$ to obtain
\begin{equation}\label{e:phi_psi}
\int_{\T \times \R} \phi_\gamma p \partial_q \psi \, \mu(dp \, dq) = 0.
\end{equation}
We introduce the decomposition
\begin{equ}[e:decompphi]
\overline{\phi}_\gamma(q) =  \int \phi_\gamma(p,q)\, \nu_\beta(dp)\;,\qquad \tilde
\phi_\gamma(p,q) = \phi_\gamma(p,q) - \overline{\phi}_\gamma(q)\;,
\end{equ}
where $\nu_{\beta}(p) = Z^{-1} \exp(- \beta p^2/2)$. Note that if $\phi$ is a function from $\R$ to $\R$
such that $\int \phi(p)\, \nu_\beta(dp) = 0$, then it follows from the spectral decomposition of the
harmonic oscillator Schr\"odinger operator that $\|\d_p \phi\|^2 \ge \beta \|\phi\|^2$, where the
norms are in $L^2(\nu_\beta)$. This inequality can be applied pointwise to $\tilde \phi_\gamma$,
so that he bound
\begin{equ}[e:Poincarephi]
\|\d_p \phi_\gamma\|^2= \|\d_p \tilde \phi_\gamma\|^2 \ge \beta \|\tilde \phi_\gamma\|^2
\end{equ}
holds.

Substituting the decomposition \eref{e:decompphi} and \eref{e:phi_psi}
into the expression for $D_\gamma$, we obtain
\begin{equs}
D_\gamma &= \int_{\T \times \R} \tilde{\phi}_\gamma p \, \mu(dp \, dq)  =  \int_{\T \times \R}
\tilde{\phi}_\gamma p (1 + \partial_q \psi) \, \mu(dp \, dq)
\\ & \leq \|\tilde{\phi}_\gamma \|\, \| p (1 + \partial_q \psi) \| \leq \sqrt{\beta^{-1}} \|\partial_p \phi_\gamma \| \|1 + \partial_q \psi \| \|p \|
\\ & \leq  \sqrt{\gamma D_\gamma} \|1 + \partial_q \psi \| \sqrt{\beta^{-1}}.
\end{equs}
Here, we used \eref{e:Poincarephi} on the second line
and we used the fact that the effective diffusion coefficient can be written as
$$
D_\gamma = \frac{1}{\gamma \beta} \|\partial_p \phi_\gamma \|^2\;.
$$
to go from the second to the third line. It follows from this calculation that $D_\gamma
\leq \frac{1}{\gamma \beta} \|1 + \d_q \psi\|^2$, so that \eqref{e:deff_a_1} follows by
taking $\psi$ in the above estimate to be $\chi$, the solution of \eqref{e:cell_smol},
and by using \eqref{e:deff_smol}.

Now we proceed with the bound from below. This time, we use a test function $\psi$ of the form
\begin{equ}
\psi(p,q) = \left\{\begin{array}{rl} \phi\circ H & \text{for $p \ge 0$} \\
-\phi \circ H& \text{for $p \le 0$} \end{array}\right.
\end{equ}
where $\phi\colon \R_+ \to \R_+$ is a smooth function with $\phi(H) = 0$ for $H
\le E_0$ and such that $\lim_{H \downarrow E_0} \phi'(H) \neq 0$. Plugging this ansatz into equation
\eqref{e:weak}, we obtain
\begin{equ}
 \beta^{-1} \int_{\T \times \R} \partial_p \phi_\gamma  \partial_p \psi \, \mu(dp \, dq) =
 \int_{\T \times \R} p \psi  \, \mu(dp \, dq) = -\beta^{-1} \int_{\T \times \R} \partial_p \psi
 \, \mu(dp \, dq).
\end{equ}
Here, we used integration by parts and the explicit expression for $\mu$ in order to obtain
the second equality.
Cauchy-Schwarz now yields:
$$
D_\gamma = {1\over \gamma\beta} \|\partial_p \phi_\gamma \|^2 \geq \frac{\left(\int_{\T \times \R} \partial_p \psi \,
 \mu(dp \, dq) \right)^2}{\gamma\beta\|\partial_p \psi \|^2}\;.
$$
At this point, we notice that, with the notations of Section~\ref{sec:critical}, this is equivalent to
\begin{equ}
D_\gamma \ge {\Bigl(2 Z_\beta^{-1} \int_{E_0}^\infty \phi'(z) e^{-\beta z}\,dz\Bigr)^2
\over 2 \gamma \beta Z_\beta^{-1} \int_{E_0}^\infty \bigl(\phi'(z)\bigr)^2 S(z) e^{-\beta z}\,dz}\;.
\end{equ}
Choosing $\phi$ such that $\phi'(z) = 1/S(z)$, we finally obtain
\begin{equ}
D_\gamma \ge {2 \over  \gamma \beta Z_\beta} \int_{E_0}^\infty {e^{-\beta z}\over S(z)}\, dz
\end{equ}
which, combined with \eref{e:deff_FW_1}, is the required bound.
\end{proof}

%
%
\subsection{Large $\pmb{\gamma}$ Asymptotic}
It is well known that, when $\gamma$ is large, solutions to the Langevin
equation~\eqref{e:lang} are approximated by solutions to the Smoluchowski
equation~\eqref{e:smoluch1}, see e.g. \cite[Thm. 10.1]{nelson}, \cite{freidlin8}. It is
therefore not surprising that a result similar to Proposition~\ref{cor:limit} holds in the large
$\gamma$ limit:
$$
\lim_{\gamma \rightarrow \infty} \gamma D = \bar D =  \frac{1}{\beta Z \widehat{Z}}\;,
$$
where
\begin{eqnarray*}
Z = \int_{0}^1 e^{-\beta V(q)} \, dq, \quad \widehat{Z} = \int_{0}^1 e^{\beta V(q)} \,
dq.
\end{eqnarray*}
In the above we used the fact that $\bar D$ can be calculated explicitly~\cite[Sec.
13.6]{PavlSt06b}. It is also quite straightforward (at least formally) to obtain the next
term in the small $\gamma^{-1}$ expansion for $D_\gamma$. We solve perturbatively \eqref{e:cell}
using the technique presented in~\cite[Ch. 8]{HorsLef84} we obtain
\begin{figure}
\begin{center}
\includegraphics[width=2.8in, height = 2.8in]{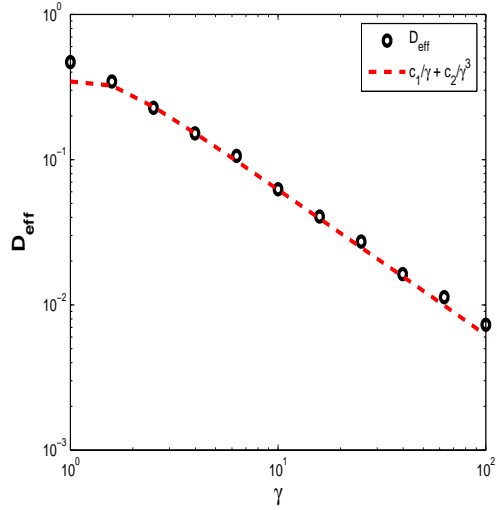}
\caption{Large $\gamma$ asymptotic for the effective diffusivity.}
\label{fig:deff_gamma}
\end{center}
\end{figure}
\begin{equation}
D_\gamma = \frac{1}{\beta \gamma Z \widehat{Z}} - \frac{\beta Z_1}{\gamma^3  Z \widehat{Z}^2} +
\mathcal{O} \left( \frac{1}{\gamma^5} \right), \label{e:deff_gamma}
\end{equation}
where $Z$ and $\hat Z$ are as before, and $Z_1$ is given by
$Z_1 = \int_0^1  (V'(q))^2 e^{\beta V(q)} \, dq$.

In Figure \ref{fig:deff_gamma}, we plot the diffusion coefficient $D$ for the nonlinear
pendulum obtained from direct numerical simulations, together with the approximation
\eqref{e:deff_gamma}. As expected, the agreement between the result of the the Monte-Carlo
simulation and the theoretical prediction is very good, even for values of $\gamma$ that are close
to $\mathcal{O}(1)$.\footnote{In the case where the potential has period $\ell$, as opposed
to $1$, then formula~\eqref{e:deff_gamma} has to be multiplied by $\ell^2$ and all the
integrals that define the coefficients that appear in the formula are taken from $0$ to
$\ell$. }
%
%
%
%
%
\section[The intermediate regime: the $\alpha \in (0,+\infty)$ case]{The intermediate regime:
the $\pmb{\alpha \in (0,+\infty)}$ case}
\label{sec:inter}

In this section we consider intermediate length and time scales, between and
Freidlin-Wentzell and the central limit theorem ones. In particular, we prove
Theorem~\ref{thm:main_2}, namely that, at intermediate length/time scales the particle
position converges weakly to a Brownian motion with the Freidlin-Wentzell effective
diffusivity~\ref{e:deff_FW_1}.

\begin{proof}[Proof of Theorem~\ref{thm:main_2}]
Notice that the stationarity assumption implies that, for every smooth function in
$L^1(\mu)$, periodic in $q$,
$$
\E f(p(t),q(t)) = \int_{\T \times \R} f(p,q) \, \mu(dp \, dq).
$$
Let $\phi_\gamma$ be the solution to the Poisson equation~\eqref{e:cell}. We apply
It\^{o}'s formula to $\phi_\gamma(p,q)$ to obtain
\begin{equs}
q^\gamma(t) & =  \lambda_\gamma q(0) + \lambda_\gamma \int_0^{t/\mu_\gamma} p(s) \, ds
\\ & =  \lambda_\gamma q(0) - \lambda_\gamma \gamma^{-1}
\Big(\phi_\gamma(p(t/\mu_\gamma), q(t/\mu_\gamma))  - \phi_\gamma(p(0),q(0))\Big) \\ & \qquad +
\sqrt{2 \gamma^{-1} \lambda_\gamma^2 \beta^{-1}} \int_0^{t/\mu_\gamma}
  \partial_p \phi_\gamma(p(s),q(s)) \, d W(s)  \\ & =:
 \lambda_\gamma q(0) +  R^\gamma + M^\gamma,
\end{equs}
where $\lambda_\gamma, \, \mu_\gamma$ are given in~\eqref{e:resc_alpha}. We obviously
have that
$
\lim_{\gamma \rightarrow 0}\E \left(\lambda_\gamma q(0) \right)^2 = 0
$.
Proposition~\ref{prop:phi}, the stationarity assumption and our assumption that $\alpha
>0$, furthermore imply that
$$
\E |R^\gamma|^2 \leq C \gamma^{\alpha} \rightarrow 0,
$$
as $\gamma \rightarrow 0$.

Consider now the martingale term $M^\gamma$. According to the martingale central limit
theorem \cite{EthKur86}, in order to prove that $M^\gamma$ converges to a Brownian
motion, it is sufficient to show that the quadratic variation process $\langle M^\gamma
\rangle$, converges weakly to a constant times $t$. This quadratic variation process is
given by:
$$
\langle M^\gamma \rangle_t = \frac{2 \lambda_\gamma^2}{ \gamma \beta}
\int_0^{t/\mu_\gamma} |\partial_p \phi_\gamma(p(s),q(s)) |^2 \, ds\;.
$$
Define
$$
f_\gamma(p,q) :=  \frac{2 \lambda_\gamma^2}{\gamma \mu_\gamma \beta} |\partial_p
\phi_{\gamma}(p,q) |^2 = 2 \beta^{-1} |\partial_p \phi_\gamma |^2
$$
and
$$\overline{f}_\gamma := \int_{\T \times \R} f_\gamma(p, q)\,\mu(dp \, dq)\;,$$
where $\mu = Z^{-1} \exp(- \beta H(p,q)) \, dp \, dq$. It follows from
Propositions~\ref{prop:phi} and~\ref{prop:resolv_norm} that $\overline{f}_\gamma$ remains
bounded between two positive constants as $\gamma \to 0$.

In order to bound the error between $\scal{M^\gamma}_t$ and $\overline{f}_\gamma t$, the
idea is to subdivide the interval $[0,t]$ into $N$ `small' intervals of size $\tau = t/N$
and to add the individual errors made at each time interval. Denote $t_k = k\tau$ and
$\eps_k = |\scal{M^\gamma}_{t_{k+1}} - \scal{M^\gamma}_{t_{k}} - \overline{f}_\gamma
\tau|$. Then, the fact that $\scal{M^\gamma}$ is an increasing process implies that
\begin{equ}[e:bounderror]
\sup_{s \in [0,t/\mu_\gamma]} |\scal{M^\gamma}_s - \overline{f}_\gamma s| \le
\sum_{k=1}^N \eps_k + \overline{f}_\gamma \tau + \sup_{k=1}^N \eps_k \le 2\sum_{k=1}^N
\eps_k + \overline{f}_\gamma \tau\;.
\end{equ}
The individual errors $\eps_k$ can be bounded in a standard way by
\begin{equs}
\E \eps_k^2 &= \E \Bigl(\int_{t_{k-1}}^{t_{k}} f_\gamma(p(s/\mu_\gamma),
q(s/\mu_{\gamma}))\, ds - \overline{f}_\gamma \tau \Bigr)^2 \\
& = \E
\Bigl(\int_{t_{k-1}}^{t_{k}} f_\gamma(p(s/\mu_\gamma), q(s/\mu_\gamma))\, ds\Bigr)^2 -
\overline{f}_\gamma^2 \tau^2 \\
& = \int_0^\tau \int_0^\tau \int_{\T\times \R} f_\gamma(p, q) \bigl(\CP_{|r-s|/\mu_\gamma}^\gamma f_\gamma\bigr)(p, q)\,\mu(dp \, dq)\,dr\,ds - \overline{f}_\gamma^2 \tau^2 \\
&\le \int_0^\tau \int_0^\tau \|f_\gamma\|^2 e^{-C \gamma |r-s|/\mu_\gamma}\,dr\,ds \le {C
\over \gamma} \|f_\gamma\|_\mu^2 \tau \mu_\gamma\;.
\end{equs}
Here, we used the stationarity assumption. We also used Theorem~\ref{theo:sgap} to bound
the action of the semigroup $\CP_t^\gamma$. Combining this with
\eref{e:bounderror} yields
\begin{equs}
\E \Bigl( \sup_{s \in [0,t/\mu_\gamma]}& |\scal{M^\gamma}_s - \overline{f}_\gamma s|
\Bigr)  \le C \|f_\gamma\|_\mu \sqrt{\mu_\gamma \over  \gamma \tau} +
\overline{f}_\gamma \tau  \\ & = C \|\partial_p \phi_\gamma \|^2_{L^4(\mu)}
\gamma^{\alpha} \tau^{-1/2} + C \|\partial_p \phi_\gamma \|^2 \tau  \leq C \left(
\gamma^{\alpha - 1/2} \tau^{-1/2} + \tau \right).
\end{equs}
We now choose $\tau = \gamma^{\zeta}$ for some $\zeta > 0$, arbitrarily small. Since
we assumed $\alpha > 1/2$, we conclude that
$$
\lim_{\gamma \to 0} \E \Bigl( \sup_{s \in [0,t/\mu_\gamma]} |\scal{M^\gamma}_s - \overline{f}_\gamma s|
\Bigr) = 0\;.
$$
Furthermore, it follows from the definition of $\bar f_\gamma$ and from Proposition~\ref{cor:limit} that
one has $\lim_{\gamma \to 0} \bar f_\gamma = 2D^*$.
This immediately implies that, as $\gamma \rightarrow 0$, $\langle M^\gamma \rangle_t$ converges
to
$2D^* t$ in $L^1(\mu)$ and therefore $M^\gamma$ converges to a
Brownian motion with diffusivity $D^*$.
\end{proof}

The proof of the result stated in Remark~\ref{rem:gamma_large} on the large $\gamma$
asymptotic is essentially identical to the one presented above: It\^{o}'s formula, our
assumption of stationarity and the scaling $\lambda_\gamma = \gamma^{-\alpha}$ lead to
$$\gamma^{-\alpha} q(t \gamma^{1 + 2 \alpha}) = M^\gamma + R^\gamma,$$
where $\lim_{\gamma \rightarrow \infty} \|R^{\gamma} \| = 0$ and
$$
M^\gamma = 2 \gamma^{-1 - 2 \alpha} \beta^{-1} \int_0^{t \gamma^{1 +2 \alpha}}
|\partial_p \phi_\gamma (p(s), q(s))|^2 \, ds.
$$
The result now follows from the martingale central limit theorem, the subdivision of
$[0,t]$ that was used in the proof above and the fact that, for $\gamma \gg 1$,
estimate~\eqref{e:estim_l4} becomes
$$
\|\partial_p \phi_\gamma \|_{L^4(\mu)} \leq C,
$$
for some constant independent of $\gamma$.
%
%
%
%
%
\section{Resolvent bounds}
\label{sec:resolvent}
In this section, we obtain the main bounds on the solution $\phi_\gamma$ of the Poisson
equation \eref{e:cell}. It will be convenient to work not only in the $L^2$ space
weighted by the invariant measure $\mu = \exp(-\beta H(p,q))\,dp\,dq$, but in the whole
scale of spaces $\CH_\delta = L^2(\mu_\delta)$ for $\delta \in (0,\beta]$. The main
technical difficulty is to obtain bounds in spaces for $\delta \neq \beta$, but this
seems to be required in order to obtain the bound on the $\L^4$-norm of $\d_p
\phi_\gamma$ required in Section~\ref{sec:inter}. We first focus on bounds for the case
$\delta = \beta$.

The norm and scalar product in $\CH_\delta$ will be
denoted by $\|\cdot\|_\delta$ and $\scal{\cdot,\cdot}_\delta$ respectively. The subscript is
omitted in the case $\delta = \beta$. We also denote by $\mu_\delta$ the probability
measure on $\T \times \R$ proportional to $\exp(-\delta H(p,q))\,dp\,dq$ and by $\nu_\delta$
the  probability measure  on $\R$ proportional to $\exp(-\delta p^2/2)\, dp$.

\subsection{Bounds in $\pmb{\CH_\beta}$}

We have the following preliminary bound:

\begin{prop}
\label{prop:phi} Let $\phi_\gamma$ denote the solution of \eqref{e:cell} and assume that
$V(q) \in C^{\infty}_{per}(\T)$. Then, $\phi_\gamma$ satisfies the bound
\begin{equation}\label{e:phi_estim1}
\| \partial_p \phi_\gamma \| \leq C\;,
\end{equation}
for some constant $C$ independent of $\gamma$.
\end{prop}
\begin{proof}
Existence and uniqueness of solutions to \eref{e:cell} is proved for example
in~\cite{papan_varadhan}, see also~\cite[Thm. 3.3]{HairPavl04}. The smoothness of the
solution follows from the hypoellipticity of the operator $L^\gamma$.
Estimate~\eqref{e:phi_estim1} follows from the Poincar\'{e} inequality for Gaussian
measures in the following way: we multiply \eqref{e:cell} by $\phi_\gamma$ and integrate by
parts on the left hand side to obtain
$$
\beta^{-1} \|\partial_p \phi_\gamma \|^2 = \scal{p, \overline{\phi} + \tilde{\phi}}\,,
$$
where we defined
\begin{equ}
\overline{\phi}_\gamma(q) =  \int \phi_\gamma(p,q)\, \nu_\beta(dp)\;,\qquad \tilde
\phi_\gamma(p,q) = \phi_\gamma(p,q) - \overline{\phi}_\gamma(q)\;.
\end{equ}
Since, for every $q$, $\tilde {\phi}_\gamma$ averages to zero with respect to
$\nu_\beta$, it satisfies the Poincar\'{e} inequality.
On the other hand, we have that
$\scal{p, \overline{\phi}_\gamma} = 0$, since $\overline{\phi}_\gamma$ is a function of $q$ only.
Hence, by Cauchy-Schwarz and Poincar\'{e}, we have the bound
\begin{equs}
\|\partial_p \phi_\gamma \|^2 & = \beta \scal{p,\tilde{\phi}_\gamma} \leq \beta \|p\|
\|\tilde{\phi}_\gamma \| \le C \| \partial_p \tilde{\phi}_\gamma \| = C \| \partial_p \phi_{\gamma}
 \|\;,
\end{equs}
for some constant $C$ independent of $\gamma$. This concludes the proof.
\end{proof}

We proceed now with the proof of Theorem~\ref{theo:sgap} which we restate here for the reader's
convenience.
\begin{theorem}\label{thm:boundSG}
There exist constants $C$ and $\alpha$ independent of $\gamma$ such that
\begin{equ}[e:sgbound2]
\|e^{L_\gamma t} f\| \le Ce^{-\alpha t} \|f\|\;,
\end{equ}
holds for every $t>0$, every $\gamma < 1$, and every  $f \in \L^2(\mu)$ such
that $\int f d \mu = 0$.
\end{theorem}
\begin{proof}
The proof is a variation on the commutator techniques introduced in
\cite{Koh69PDO,Koh77}
and further developed in
\cite{EPR99, EH00,HerNie02IHA,HelNie05HES,Vil04HPI}. The argument given here is
actually mainly inspired by the techniques developed by Villani in \cite{Vil04HPI}.
 In particular, we make use of his idea of constructing a `skewed' scalar product in which
 the coercivity of $L_\gamma$ becomes apparent. The main
difference is that we are going to track carefully the dependence of the various
terms on the parameter $\gamma$.

We will use the following notations:
\begin{equs}[2]
A &= \beta^{-1/2}\d_p \;,\quad&\quad A^* &= -\beta^{-1/2} \d_p + \beta^{1/2}p\;,\\
B &= p\d_q - V'(q)\d_p\;, \quad&\quad B^* &= -B\;,
\end{equs}
The reason why we are using the symbol $B$ for the Liouville operator and not $\cA$ as previously
is to be consistent with the notations adopted in \cite{Vil04HPI}.
With these notations, we have $L_\gamma = - A^* A + \gamma^{-1} B$.
Here, the adjoints $A^*$ and $B^*$ are taken with respect to the scalar product in $\CH = \L^2(\mu)$.
We also introduce the operators
\begin{equs}
\com  &= [A,B] = \beta^{-1/2}\d_q\;,\\
R &= [\com ,B] = - \beta^{-1/2}V''(q)\d_p = -V''(q) A\;.
\end{equs}

We introduce the symmetric sesquilinear form $\sscal{\cdot,\cdot}$ defined by polarisation
from
\begin{equ}
\sscal{f,f} = a\scal{f,f} + \gamma \bigl(b\scal{Af, Af} + 2\,\Re \scal{Af, \com f} +
b\scal{\com f, \com f}\bigr)\;,
\end{equ}
for some constants $a$ and $b$ to be determined later. If we take $b > 1$, then this
is indeed positive definite and induces a norm equivalent to the norm $\|\cdot\|_{1,\gamma}$ given by
\begin{equ}[e:defnorm]
\|f\|_{1,\gamma}^2 = \|f\|^2 + \gamma \bigl( \|\d_p f\|^2 +  \|\d_q f\|^2\bigr)\;.
\end{equ}
Following the same manipulations as in  \cite[Theorem~18]{Vil04HPI},
we see that there exists a constant $c$ independent of $\gamma$ such that
\begin{equs}
\Re \scal{f, L_\gamma f} &= -\|Af\|^2\;, \\
\Re \scal{Af, AL_\gamma f} &\le -\|A^2 f\|^2 + c \|Af\|^2 + {c \over \gamma} \|A f\| \|\com  f\|\;,\\
\Re \scal{Af, \com L_\gamma f} + \Re \scal{AL_\gamma f, \com f} &\le -{1\over \gamma} \|\com  f\|^2 + {c\over \gamma} \|Af\|^2 \\
&\qquad + {c} \|A^2 f\| \|\com A f\| + c \|\com f\| \|Af\|\;,\\
\Re \scal{\com f, \com L_\gamma f} &\le - \|\com A f\|^2 + {c\over \gamma} \|A f\| \|\com  f\|\;.
\end{equs}
It is now easy to see that we can choose $a \gg b \gg 1$ sufficiently large (but still independently of $\gamma$!) so that
\begin{equ}[e:boundL]
\Re \sscal{f,L_\gamma f} \le - \|Af\|^2 - \|\com f\|^2 - \gamma \bigl(\|A^2 f\|^2 + \|\com Af\|^2\bigr)\;.
\end{equ}
Note now that, provided that $f$ is centred with respect to $\mu$, the Poincar\'e inequality
tells us that there exists a constant $\kappa$
such that
\begin{equ}
 \|Af\|^2 + \|\com f\|^2 \ge \kappa \|f\|^2\;,
\end{equ}
so that \eref{e:boundL} implies in particular that $\Re \sscal{f,L_\gamma f} \ge \kappa' \sscal{f,f}$
for some $\kappa'$. This immediately implies that there exist positive constants $C$ and $\alpha$
such that
\begin{equ}[e:gapH1]
\|e^{L_\gamma t} f\|_{1,\gamma} \le Ce^{-\alpha t} \|f\|_{1,\gamma}\;.
\end{equ}
We will now show that there exists a time $\tau$ and a constant $C$, both independent of $\gamma$
(provided that $\gamma$ is sufficiently small) such that
\begin{equ}[e:regul]
\|e^{L_\gamma \tau} f\|_{1,\gamma} \le C \|f\|\;.
\end{equ}
Combining this with \eref{e:gapH1} and \eref{e:defnorm} then implies that \eref{e:sgbound2} holds.

In order to show \eref{e:regul}, we combine the previous technique
with the usual trick for proving regularisation results for parabolic PDEs.
We fix a smooth function $f$ and we define the quantity
\begin{equ}
2 A_f(t) = K \|f_t\|^2 + \gamma (t \|A f_t\|^2 + t^3\|\com  f_t\|^2 + \delta t^2 \scal{A f_t, \com  f_t})\;,
\qquad f_t = \CP_t^\gamma f \;,
\end{equ}
for some (large) constant $K$ and some (small) constant $\delta$
to be determined later. Taking the time derivative of $A_f$, we obtain
\begin{equs}
\d_t A_f &\le - K \|A f_t\|^2 + \gamma \|A f_t\|^2 + t \bigl(- \gamma \|A^2 f_t\|^2 + c\gamma\|A f_t\|^2 + c\|A f_t\| \|\com  f_t\|\bigr) \\
&\quad + 3\gamma t^2 \|\com  f_t\|^2 + t^3 \bigl(- \gamma \|\com A f_t\|^2 + c \|A f_t\| \|\com  f_t\| \bigr)
+ 2 t \gamma \delta \scal{A f_t, \com f_t} \\
& \quad  +\delta  t^2 \bigl(- \|\com f_t \|^2 + c \|Af_t\|^2 + c \gamma \|A^2 f_t\| \|\com A f_t\| + c \gamma \|\com f_t\| \|Af_t\|\bigr) \\
& \le - {K\over 2} \|A f_t\|^2  - \gamma t \|A^2 f_t\|^2 - \gamma t^3 \|\com A f_t\|^2 - \delta t^2 \|\com  f_t\|^2 \\
&\quad + ct \|A f_t\| \|\com  f_t\| + c\gamma \delta t^2 \|A^2 f_t\| \|\com A f_t\|\;.
\end{equs}
(For the second inequality we changed the value of the constant $c$ and we assumed
that $t \in [0,1]$.)
Notice now that one can first choose $\delta$ sufficiently small (but independently of $\gamma$)
such that
\begin{equ}
 c\gamma \delta t^2 \|A^2 f_t\| \|\com A f_t\| \le  \gamma t \|A^2 f_t\|^2 + \gamma t^3 \|\com A f_t\|^2\;.
\end{equ}
We can then choose $K$ sufficiently large so that
\begin{equ}
ct \|A f_t\| \|\com  f_t\| \le {K\over 2} \|A f_t\|^2 + \delta t^2 \|\com  f_t\|^2 \;.
\end{equ}
With these choices, we get $\d_t A_f \le 0$, so that $A_f(1) \le A_f(0)$. This immediately implies
the bound \eref{e:regul}.
\end{proof}

By simply integrating from $0$ to $\infty$, this implies that one has the resolvent bound:
\begin{equ}[e:resolventbound]
\|L_\gamma^{-1} f\| \le C \|f\|\;,
\end{equ}
holding for every $f \in \CH_\beta$ such that $\scal{1, f} = 0$. It turns out that, up to a constant,
this bound is actually optimal:

\begin{prop}
\label{prop:resolv_norm}
There exists a constant $C$ independent of $\gamma$ such that the operator norm of the resolvent satisfies
$$
\|L_\gamma^{-1}\| \ge C\;.
$$
\end{prop}

\begin{proof}
We make use of the fact that the norm of the resolvent can be characterized by
\begin{equation}
\|L_\gamma^{-1} \| = \left( \inf_{f \in \CD(L_\gamma)\,:\, \scal{1,f} = 0}
\frac{\|L_\gamma f \|}{\|f \|} \right)^{-1}\;.
\label{e:res_norm}
\end{equation}
If we take $f$ of the form $f = \phi \circ H$ for an arbitrary smooth bounded function $\phi$
such that $\scal{1,f} = 0$, then $L_\gamma f$ (and therefore also $\|L_\gamma f\|$) is
independent of $\gamma$. Similarly, $\|f\|$ is independent of $\gamma$ so that the infimum
appearing in \eref{e:res_norm} is bounded from above by a constant independent of $\gamma$,
thus proving the claim.
\end{proof}

We now use these estimates to obtain bounds on the solution $\phi_{\gamma}$ to the
Poisson equation $-L_\gamma \phi_{\gamma} = p$.

\begin{proposition}\label{prop:goodbounds}
There exists a constant $C$ independent of $\gamma$ such that
\begin{equ}
\|\phi_{\gamma}\| + \|\d_p \phi_{\gamma}\| + \|\d_q \phi_{\gamma}\|  \le C\;.
\end{equ}
\end{proposition}

\begin{proof}
We have from \eref{e:boundL}  and the fact that the $\|\cdot\|_{1,\gamma}$-norm of
$p$ is bounded independently of $\gamma$ that
\begin{equ}
\|\d_p\phi_{\gamma}\|^2 + \|\d_q\phi_{\gamma}\|^2 \le \Re\sscal{\phi_{\gamma},
L\phi_{\gamma}} = \Re\sscal{\phi_{\gamma}, p} \le C \|\phi_{\gamma}\|_{1,\gamma} \le C\;.
\end{equ}
The last inequality followed from the bound \eref{e:gapH1}.
One can actually extract slightly more from the above bounds.
At this stage, we note that one can write $B = \gamma(L - A^*A)$, so that
\begin{equ}
\|B\phi_{\gamma}\|^2 = \gamma \scal{p,B\phi_{\gamma}} - \gamma \scal{A\phi_{\gamma},
AB\phi_{\gamma}}\;.
\end{equ}
Since furthermore $B$ is antisymmetric and $[A,B] = \d_q$, we have
\begin{equ}
\|B \phi_{\gamma}\|^2 \le C\gamma + C\gamma \|\d_p \phi_{\gamma}\| \|\d_q
\phi_{\gamma}\|\;.
\end{equ}
Collecting this with the previous estimates, we obtain the existence of a constant $C$
such that $\|B \phi_{\gamma}\| \le C\sqrt \gamma$, which in turn yields a bound of the
type $\|A^* A \phi_{\gamma}\| \le C/\sqrt{\gamma}$.
\end{proof}
%
%
\subsection{Bounds in $\pmb{\CH_\delta}$ with $\pmb{\delta \neq \beta}$}

We now show that similar bounds hold in every $\CH_\delta$. The main difficulty is that these
spaces are no longer weighted by the invariant measure of the system, so that several
simplifications are lost. In particular, the very useful relation
$\Re \scal{\phi, L_\gamma \phi} = -\beta^{-1}\|\d_p \phi\|^2$ does not hold anymore.

Throughout this section, we will write $L_\sym$ for the symmetric part of $L_\gamma$ in
$\CH_\delta$: $\scal{\phi,L_\gamma \phi}_\delta = \scal{\phi, L_\sym \phi}_\delta$. An
explicit calculation shows that one has
\begin{equ}
L_\sym = - \beta^{-1} \d_p^* \d_p + {\beta - \delta \over 2\beta} - {\delta
(\beta-\delta)\over 2\beta} p^2\;,
\end{equ}
where we denote by $\d_p^* = -\d_p + \delta p$ the adjoint of $\d_p$ in $\CH_\delta$.
Note that one has $L_\sym = - \beta^{-1} \d_p^* \d_p$ if and only if $\delta = \beta$. A
standard calculation shows that $L_\sym$ is unitarily equivalent to the Schr\"odinger
operator corresponding to the harmonic oscillator, so that one can explicitly compute its
spectral decomposition. In order to do so, we define
\begin{equ}[e:defA]
\alpha = -{\delta \over 2} + {\sqrt{\delta (2\beta - \delta)} \over 2}\;,\qquad A =
\beta^{-1/2}\bigl(\d_p + \alpha p\bigr)\;,
\end{equ}
and we note that one can write
\begin{equ}[e:Lsym]
L_\sym = - A^* A + {\beta - \sqrt{\delta (2\beta - \delta)} \over 2\beta}\;.
\end{equ}
Furthermore, one has $[A^*, A] = -(2\alpha + \delta)/\beta$. This shows that the
eigenvalues of $L_\sym$ are given by
\begin{equ}
\lambda_n = {\beta - \sqrt{\delta (2\beta - \delta)} \over 2\beta} -
n {\sqrt{\delta (2\beta - \delta)} \over \beta}\;,\qquad n=0,1,\ldots\;,
\end{equ}
and the corresponding eigenfunctions are $f_n \propto (A^*)^n f_0$ with $f_0 \propto
\exp(-\alpha p^2/2)$. In the special case $\delta = \beta$, one simply has $\lambda_n = -n$.

In order to obtain bounds in $\CH_\delta$, we note that \eref{e:cell} yields the relation
\begin{equ}
 \|A \phi_{\gamma}\|_\delta^2 - {\beta - \sqrt{\delta (2\beta - \delta)} \over 2\beta}
 \|\phi_{\gamma}\|_\delta^2 = \scal{\phi_{\gamma},p}_\delta\;.
\end{equ}
This immediately implies that there exist constants $C$ and $N$ independent of $\gamma$ such that
one has
\begin{equ}
{1\over 2}\bigl(\|A \phi_{\gamma}\|_\delta^2 + \|\phi_{\gamma}\|_\delta^2\bigr) \le
|\scal{\phi_{\gamma} ,  p}_\delta| + C \sum_{k=0}^N  \int_\T \Bigl(\int_\R f_k(p)
\phi_{\gamma}(p,q)e^{-\delta p^2/2}\, dp\Bigr)^2\,dq\;.
\end{equ}
On the other hand, all the $f_k$'s are of the form $P_k(p) e^{-\alpha p^2/p}$ for some polynomial
$P_k$ of degree $k$. This implies that, for every $\delta' < 2(\alpha + \delta)$, there exist
constants $C_1$, $C_2$ such that one has the bound
\begin{equ}
\|A \phi_{\gamma}\|_\delta + \|\phi_{\gamma}\|_\delta \le C_1 + C_2
\|\phi_{\gamma}\|_{\delta'}\;.
\end{equ}
Since, for $\delta \ge \beta$, one has $\alpha \ge 0$, one has in particular the bound
\begin{equ}[e:boundAPhi]
\|A \phi_{\gamma}\|_\delta + \|\phi_{\gamma}\|_\delta \le C \bigl(1 +
\|\phi_{\gamma}\|_{2\delta}\bigr)\;.
\end{equ}
This calculation shows that:

\begin{proposition}\label{prop:boundPhiHdelta}
One has $\phi_{\gamma} \in \bigcap_{\delta > 0} \CH_\delta$ and, for every $\delta \in
(0,\beta]$, there exists a constant $C$ independent of $\gamma$ such that
\begin{equ}[e:boundPhiHdelta]
\|\d_p \phi_{\gamma}\|_\delta^2 
+ \|\phi_{\gamma}\|_\delta^2 \le C\;.
\end{equ}
\end{proposition}

\begin{proof}
Since one has $\|\phi_{\gamma}\|_\delta \le C\|\phi_{\gamma}\|_{\delta'}$ for $\delta \ge
\delta'$, we can apply \eref{e:boundAPhi} recursively to obtain
\begin{equ}
\|A \phi_{\gamma}\|_\delta + \|\phi_{\gamma}\|_\delta \le C \bigl(1 +
\|\phi_{\gamma}\|_{\beta}\bigr) \le C\;,
\end{equ}
where we made use of Proposition~\ref{prop:phi} for the second inequality (the two
constants $C$ are of course different). Since furthermore $\|p \phi_{\gamma}\|_\delta \le
C\|\phi_{\gamma}\|_{\delta'}$ for $\delta \ge \delta'$, this proves the claim.
\end{proof}

We are now going to show that it is also possible to obtain an order $1$ bound for
$\|\d_q \phi_{\gamma}\|_\delta$, but as before this is less straightforward. We first
start with the following preparatory result:
\begin{proposition}\label{prop:secondbounds}
There exists a constant $C$ independent of $\gamma$ such that
\begin{equ}
\|\d_p^2 \phi_{\gamma}\| + \|\d_p\d_q \phi_{\gamma}\| + \|\d_q^2 \phi_{\gamma}\| \le {C
\over \sqrt \gamma}\;.
\end{equ}
\end{proposition}

\begin{proof}
The bound for $\|\d_p^2 \phi_{\gamma}\|$ was obtained in
Proposition~\ref{prop:goodbounds}. In order to obtain the bound on $\d_p\d_q
\phi_{\gamma}$, note that one has $[L_\gamma, \d_q] = \gamma^{-1} V''(q) \d_p$, so that
$L_\gamma \d_q \phi_{\gamma} = \gamma^{-1} V''(q) \d_p \phi_{\gamma}$. Therefore, we have
\begin{equ}
\|\d_p\d_q \phi_{\gamma}\|^2 = \scal{\d_q \phi_{\gamma}, L_\gamma \d_q \phi_{\gamma}} =
\gamma^{-1}\scal{\d_q \phi_{\gamma}, V''(q) \d_p \phi_{\gamma}} \le {C\over \gamma}
\|\d_q \phi_{\gamma}\|\|\d_p \phi_{\gamma}\|\;,
\end{equ}
so that the bound follows from Proposition~\ref{prop:goodbounds}. Finally, it follows from
\eref{e:boundL} that we have the bound
\begin{equs}
\|\d_q^2 \phi_{\gamma}\|^2 &\le |\sscal{\d_q\phi_{\gamma}, L_\gamma \d_q \phi_{\gamma}}| = \gamma^{-1} |\sscal{\d_q \phi_{\gamma}, V''(q) \d_p \phi_{\gamma}}| \\
&\le {C\over \gamma} \|\d_q \phi_{\gamma}\|\|\d_p \phi_{\gamma}\| + C \bigl(\|\d_p \d_q \phi_{\gamma}\|^2 + \|\d_p\d_q\phi_{\gamma}\|\|\d_p^2\phi_{\gamma}\| + \|\d_q^2 \phi_{\gamma}\|\|\d_p\d_q\phi_{\gamma}\| \\
&\qquad + \|\d_q^2 \phi_{\gamma}\| \|\d_p^2 \phi_{\gamma}\| + \|\d_q^2 \phi_{\gamma}\| \|\d_p \phi_{\gamma}\| + \|\d_p\d_q \phi_{\gamma}\| \|\d_p \phi_{\gamma}\|\bigr) \\
&\le {C\over \gamma} + {1\over 2}\|\d_q^2 \phi_{\gamma}\|^2\;,
\end{equs}
where we made use of all of the the previously obtained bounds in the last step.
\end{proof}

Our aim now is to mimic the proof of Theorem~\ref{theo:sgap}, with the space
$\CH_\beta$ replaced by $\CH_\delta$ for some arbitrary $\delta \in (0,\beta]$.
We define $A$ as in \eref{e:defA} and we set as before $B = p\d_q - V'(q)\d_p$.
We furthermore define the
operator $\tilde B$ (which is antisymmetric in $\CH_\delta$) by
\begin{equ}
\tilde B = {\delta -\beta \over \beta}\Bigl(p\d_p - {\delta \over 2} p^2 + {1\over 2}\Bigr)\;.
\end{equ}
With these notations, we can check that $L_\gamma$ can be written as
\begin{equ}
L_\gamma = - A^* A + {1\over \gamma} B + \tilde B + {\beta - \sqrt{\delta(2\beta - \delta)}\over 2\beta}\;,
\end{equ}
where the adjoint of $A$ is taken in $\CH_\delta$. This motivates the definition of an operator
$\hat L_\gamma$ given by
\begin{equ}
\hat L_\gamma = - A^* A + {1\over \gamma} B \;.
\end{equ}
Our strategy is then to obtain a bound similar to \eref{e:boundL} with $L_\gamma$ replaced by $\hat L_\gamma$ and to use make use of the fact that the difference between $L_\gamma$ and
$\hat L_\gamma$ is sufficiently ``small''.

\begin{theorem}\label{theo:ultimatebounds}
For every $\delta \in (0,\beta]$, there exists a constant $C$ such that
\begin{equ}[e:boundPhidelta]
\|\d_p \phi_{\gamma}\|_\delta^2 + \|\d_q \phi_{\gamma}\|_\delta^2 + \gamma \bigl(\|\d_p^2
\phi_{\gamma}\|_\delta^2 + \|\d_p \d_q \phi_{\gamma}\|_\delta^2 + \|\d_q^2
\phi_{\gamma}\|_\delta^2\bigr) \le C\;,
\end{equ}
independently of $\gamma$.
\end{theorem}

\begin{remark}
In terms of $L^2$-estimates, these bounds are likely not to be absolutely optimal. In the
limit $\gamma \to 0$ the solution $\phi_{\gamma}$ of the Poisson equation \eref{e:cell}
indeed converges to a function of the form
\begin{equ}
\phi_{\gamma,0}(p,q) = \left\{\begin{array}{rl} \phi\circ H & \text{for $p \ge 0$} \\
-\phi \circ H& \text{for $p \le 0$} \end{array}\right.
\end{equ}
where $\phi\colon \R_+ \to \R_+$ is a smooth function with $\phi_{\gamma}(H) = 0$ for $H
\le E_0$ and such that $\lim_{H \downarrow E_0} \phi'(H) \neq 0$. Note that the second
derivative of $\phi_{\gamma,0}$ is therefore not square-integrable. For small values of
$\gamma$, it is believed \cite{Sowers03, Sowers05} that, around $H = E_0$, the function
$\phi_{\gamma}$ develops a `boundary layer' of width $\sqrt \gamma$ on which, from simple
scaling arguments, its second derivative should be of order $\gamma^{-1/2}$ .
\end{remark}

\begin{proof}[Proof of Theorem~\ref{theo:ultimatebounds}]
We define as before the operators $\com$ and $R$ by
\begin{equ}
\com = [A,B] = \beta^{-1/2} \bigl(\d_q + \alpha V'(q)\bigr) \;,
\qquad R = [\com,B] = - V''(q) A\;.
\end{equ}
With these notations, we define similarly as before the scalar product
\begin{equ}
\sscal{f,f}_\delta = a\|f\|_\delta^2 + \gamma \bigl(b \scal{Af, Af}_\delta +
2\Re \scal{Af, \com f}_\delta + b\scal{\com f,\com f}_\delta\bigr)\;,
\end{equ}
where $a$ and $b$ are constants to be determined later. Since the algebraic relations
between $\hat L_\gamma$, $A$, $B$, $\com $, and $R$ are exactly the same as above, we can
retrace step by step the proof of \eref{e:boundL} to get
\begin{equ}[e:boundLdelta]
\|A f\|_\delta^2 + \|\com  f\|_\delta^2 + \gamma \bigl(\|A^2 f\|_\delta^2 + \|\com A
f\|_\delta^2\bigr) \le - \Re \sscal{f, \hat L_\gamma f}_\delta\;.
\end{equ}
Since furthermore we know from Proposition~\ref{prop:boundPhiHdelta} that
$\|\phi_{\gamma}\|_\delta$ and (and therefore also $\|p^n\phi_{\gamma}\|_\delta$ for
every $n$) are bounded by constants independent of $\gamma$, this implies the existence
of a constant $C$ such that
\begin{equs}
\|\d_p \phi_{\gamma}\|_\delta^2 + \|\d_q \phi_{\gamma}\|_\delta^2 &+ \gamma \bigl(\|\d_p^2
\phi_{\gamma}\|_\delta^2 + \|\d_p \d_q \phi_{\gamma}\|_\delta^2\bigr) \\
&\le C\bigl(1 +
\sscal{\phi_{\gamma},\phi_{\gamma}}_\delta + |\sscal{\phi_{\gamma}, p\d_p
\phi_{\gamma}}_\delta|\bigr)\;.
\end{equs}
However, it is a straightforward calculation to check that, from
Proposition~\ref{prop:boundPhiHdelta} and the definition of $\sscal{\cdot,\cdot}_\delta$,
one has
\begin{equ}
|\sscal{\phi_{\gamma}, p\d_p \phi_{\gamma}}_\delta| \le {1\over 2}
\bigl(\|\d_p\phi_{\gamma}\|_\delta^2 + \gamma \bigl(\|\d_p^2 \phi_{\gamma}\|_\delta^2 +
\|\d_p \d_q \phi_{\gamma}\|_\delta^2\bigr)\bigr) + C\bigl(1+ \sscal{\phi_{\gamma},
\phi_{\gamma}}\bigr)\;,
\end{equ}
so that we get the bound
\begin{equ}
\|\d_p \phi_{\gamma}\|_\delta^2 + \|\d_q \phi_{\gamma}\|_\delta^2 + \gamma \bigl(\|\d_p^2
\phi_{\gamma}\|_\delta^2 + \|\d_p \d_q \phi_{\gamma}\|_\delta^2\bigr) \le C\bigl(1 +
\sscal{\phi_{\gamma},\phi_{\gamma}}_\delta\bigr)\;.
\end{equ}
We can actually even get slightly better than that, in the same way as in
Proposition~\ref{prop:secondbounds}. Using \eref{e:boundLdelta} and the commutation
relation $[\hat L_\gamma, \d_q] = \gamma^{-1} V''(q) \d_p$, we have:
\begin{equ}
\|\d_q^2 \phi_{\gamma}\|_\delta^2 \le C +  |\sscal{\d_q \phi_{\gamma}, \hat L_\gamma \d_q
\phi_{\gamma}}_\delta| \le C + C |\gamma^{-1}\sscal{\d_q \phi_{\gamma}, V''(q) \d_p
\phi_{\gamma}}_\delta|\;,
\end{equ}
so that we finally get the existence of a constant $C$ such that
\begin{equ}[e:boundPhideltaalmost]
\|\d_p \phi_{\gamma}\|_\delta^2 + \|\d_q \phi_{\gamma}\|_\delta^2 + \gamma \bigl(\|\d_p^2
\phi_{\gamma}\|_\delta^2 + \|\d_p \d_q \phi_{\gamma}\|_\delta^2 + \|\d_q^2
\phi_{\gamma}\|_\delta^2\bigr) \le C\bigl(1 +
\sscal{\phi_{\gamma},\phi_{\gamma}}_\delta\bigr)\;.
\end{equ}

Our aim now is to show that, for every $\delta \in (0,\beta]$, there exists a constant
$C$ such that $\sscal{\phi_{\gamma},\phi_{\gamma}}_\delta \le C$, independently of
$\gamma$, which will then conclude the proof of the theorem. This will be performed
thanks to a bootstrapping argument similar to the one we used already in the proof of
Proposition~\ref{prop:boundPhiHdelta}. One has, for some constant $c>0$,
\begin{eqnarray*}
c \sscal{\phi_{\gamma},\phi_{\gamma}}_\delta &\le & \|\phi_{\gamma}\|_\delta^2 +
\gamma \bigl(\|\d_p \phi_{\gamma}\|_\delta^2 + \|\d_q \phi_{\gamma}\|_\delta^2\bigr) \\
& \le & \|\phi_{\gamma}\|_\delta^2 + \gamma \int_\T \int_\R \bigl((\d_p \phi_{\gamma})^2
+ (\d_q\phi_{\gamma})^2\bigr) \,\mu_\delta ( dp\,dq ) \\
&\le & \|\phi_{\gamma}\|_\delta^2 + \gamma \int_\T \int_\R \Big( |\phi_{\gamma} \d_p^2
\phi_{\gamma}| + |\phi_{\gamma} \d_q^2 \phi_{\gamma}| \\ && + \delta |p\phi_{\gamma} \d_p
\phi_{\gamma}| + \delta |V'(q) \phi_{\gamma} \d_q\phi_{\gamma}| \Big)\, \mu_\delta (
dp\,dq )
\end{eqnarray*}
Writing $2\delta = \delta_1 + \delta_2$ with $\delta_i > 0$ and applying Cauchy-Schwarz,
we obtain from this the existence of positive constants $c$ and $C$ such that
\begin{equs}
c\sscal{\phi_{\gamma},\phi_{\gamma}}_\delta &\le\|\phi_{\gamma}\|_\delta^2 + C
\gamma\bigl(\|\phi_{\gamma}\|_{\delta_1} + \|p\phi_{\gamma}\|_{\delta_1}\bigr) \\
&\qquad \times
\bigl(\|\d_p^2 \phi_{\gamma}\|_{\delta_2} + \|\d_q^2 \phi_{\gamma}\|_{\delta_2} + \|\d_p
\phi_{\gamma}\|_{\delta_2} + \|\d_q \phi_{\gamma}\|_{\delta_2}\bigr)\;.
\end{equs}
It therefore follows from \eref{e:boundPhideltaalmost} and Proposition~\ref{prop:boundPhiHdelta} that
there is a constant (depending on the choice of $\delta_i$) such that
\begin{equ}
\sscal{\phi_{\gamma},\phi_{\gamma}}_\delta \le C \bigl(1 +
\sscal{\phi_{\gamma},\phi_{\gamma}}_{\delta_2} \bigr)\;.
\end{equ}
Since $\delta_2$ can be chosen larger than $\delta$ (actually up to, but not including
$2\delta$), we can apply this inequality recursively to bound
$\sscal{\phi_{\gamma},\phi_{\gamma}}_\delta$ by
$\sscal{\phi_{\gamma},\phi_{\gamma}}_\beta$ which in turn has already been bounded in
Propositions~\ref{prop:boundPhiHdelta} and \ref{prop:secondbounds}.
\end{proof}

As a simple corollary, we have the following bound on $\d_p \phi_{\gamma}$ which is used
in the proof of  Theorem~\ref{thm:main_2}:

\begin{corollary}
The function $\d_p \phi_{\gamma}$ belongs to $L^4(\mu_\delta)$ for every $\delta \in
(0,\beta]$ and every $\gamma>0$. Its norm in $L^4(\mu_\delta)$ is of order
$\CO(\gamma^{-1/4})$.
\end{corollary}

\begin{proof}
Denoting by $\Delta$ the Laplacian on $\R^2$, it follows from the fractional Sobolev
inequalities that
\begin{equs}
\Bigl(\int \bigl|\d_p \phi_{\gamma}\bigr|^4 &\, \mu_\delta (dp\,dq ) \Bigr)^{1/2}
\le \int \Bigl((1- \Delta)^{1/4} \d_p \phi_{\gamma} \, e^{-\delta H(p,q)/4} \, dp \, dq
\Bigr)^2 \, dp\,dq \\
&\le \int \Bigl((1- \Delta)^{1/2} \d_p \phi_{\gamma} e^{-\delta H(p,q)/4}\Bigr)\d_p
\phi_{\gamma} e^{-\delta H(p,q)/4} \, dp\,dq \\
&\le C\|\d_p \phi_{\gamma}\|_{\delta/2} \bigl(\|\d_p \phi_{\gamma}\|_{\delta/2} +
\|\d_p^2
\phi_{\gamma}\|_{\delta/2} + \|\d_q \d_p \phi_{\gamma}\|_{\delta/2} + \|p\d_p \phi_{\gamma}\|_{\delta/2}\bigr) \\
& \le C \left(1 + \gamma^{-1/2} \right) \; ,
\end{equs}
were we used Theorem~\ref{theo:ultimatebounds}.
\end{proof}

\begin{remark}
In a similar way, one can obtain, for every $p \in [1,\infty)$, bounds of order $\CO(1)$
for $\phi_{\gamma}$ $L^p(\mu_\delta)$. Unfortunately, using the Sobolev bounds obtained
for $\phi_{\gamma}$ in this section, it is not possible to obtain bounds of order
$\CO(1)$ for $\d_p \phi_{\gamma}$ in $L^p(\mu_\delta)$, even though
 we conjecture that such bounds hold true.
\end{remark}

\bibliography{./mybib}
\bibliographystyle{./Martin}
\end{document}